\documentclass[12pt]{article}
\usepackage{amsmath,amsfonts,times} 

\date  {{\small Submitted Aug. 29, 1998; 
revised March 30, 1999.}}

\input epsf     
\def\be{\begin{equation}}
\def\ee{\end{equation}}
\def\bea{\begin{eqnarray}}
\def\eea{\end{eqnarray}}

\newcommand{\eq}[1]{eq.~(\ref{#1})}  
\newcommand{\lamst}{{\lambda^*}}

\def\abs#1{\left\vert #1 \right\vert}
\def\blackbox{{\vrule height 1.3ex width 1.0ex depth -.2ex} \hskip 1.5truecm}
\def\C{{\mathcal C}}   
\def\c{\mathrm{const.}}

\def\deg  {{\rm deg}}
\def\diam{\operatorname{diam}}   
\def\dim{\operatorname{dim}}  
\def\dist  {{\rm dist}}
\def\E{{\mathrm E}}    
\def\eps {\varepsilon}
\def\F{{\mathcal F}}   

\def\limsup{\mathop{\overline{\rm lim}}} 

\def\Om{\Omega}    

\def\P{\operatorname{Prob}} 
\def\R{{\mathbb R}}    
\def\S{{\mathcal S}}   
\def\sp{\dot{\mathbb R}}   
\def\T{{\mathcal T}}   
\def\too#1{\parbox[t]{.4in} {$\longrightarrow\\[-9pt] {\scriptstyle #1}$}}
\def\V{{\mathcal V}}   
\def\Z{{\mathbb Z}}    

\newcounter{masectionnumber}
\newcommand{\masect}[1]{\setcounter{equation}{0}
\refstepcounter{masectionnumber} \vspace{1truecm plus 1cm} \noindent
    {\large\bf \arabic{masectionnumber}. #1}\par \vspace{.2cm}
      \addcontentsline{toc}{section}{\arabic{masectionnumber}. #1}
    }
 \renewcommand{\theequation}
    {\mbox{\arabic{masectionnumber}.\arabic{equation}}}

    \newcounter{masubsectionnumber}[masectionnumber]
\newcommand{\masubsect}[1]{
    \refstepcounter{masubsectionnumber} \vspace{.5cm} \noindent
  {\large\em \arabic{masectionnumber}.\alph{masubsectionnumber} #1}
    \par\vspace*{.2truecm}

\addcontentsline{toc}{subsection}
{\arabic{masectionnumber}.\alph{masubsectionnumber}\hspace{.1cm} #1}
    }

\newtheorem{lem}{Lemma}[masectionnumber]
\newtheorem{thm}[lem]{Theorem}

\newtheorem{cor}[lem]{Corollary}
\newtheorem{df}{Definition}[masectionnumber]

\newenvironment{proof_claim}{
    \addtolength{\leftmargini}{-1em}\begin{quotation} \noindent
    {\bf Proof of Claim:}}{
    \hspace*{\fill} $\diamond$
    \end{quotation}\addtolength{\leftmargini}{1em}}
\newenvironment{proof_of}[1]{
    \noindent{\bf Proof of #1:} \hspace*{1em} }{
    \hfill \blackbox\bigskip}
\newenvironment{proof}{\noindent{\bf Proof:}
   \hspace*{1em}}{\hfill \blackbox\bigskip}

\def\remark{ \noindent {\bf Remark}\hspace*{1em}  }

\newcommand{\startappendix}{ \setcounter{masectionnumber}{0}
 \renewcommand{\theequation}
    {\mbox{\Alph{masectionnumber}.\arabic{equation}}}
 \renewcommand{\thelem}
    {\mbox{\Alph{masectionnumber}.\arabic{lem}}}
\addcontentsline{toc}{section}{Appendix }
  }
\newcommand{\maappendix}[1]{          
    \setcounter{equation}{0}
\refstepcounter{masectionnumber} \vspace{1truecm plus 1cm} \noindent
    {\large\bf \Alph{masectionnumber}. #1}\par \vspace{.2cm}

    \addcontentsline{toc}{section}{\Alph{masectionnumber}. #1 }   }

    \newcounter{masubapp}[masectionnumber]
\begin{document}

\title{\vspace*{-.35in}
Scaling Limits for Minimal and Random \\
        Spanning Trees in Two Dimensions }

\author{Michael Aizenman\thanks{
Departments of Physics and Mathematics,
    Princeton University; Jadwin Hall, Princeton, NJ 08544.}
\ \ \ \
Almut Burchard\thanks{University of Virginia, Department of Mathematics,
Kerchof Hall, Charlottesville, VA 22903.
Research conducted while at Department of Mathematics, Fine Hall,
Princeton, NJ 08544.}
\\
Charles M. Newman\thanks{\noindent Courant
             Institute of Mathematical Sciences, New York
            University, 251 Mercer St., New York, NY 10012.}
\ \ \ \
David B. Wilson\thanks{
Microsoft, One Microsoft Way 31/3348, Redmond, WA 98052.  Research
conducted in large part while at the Institute for Advanced Study,
Princeton.   }
   }
\maketitle
\thispagestyle{empty}        

\begin{abstract} A general formulation is presented for
continuum scaling limits of stochastic spanning trees.
A spanning tree is expressed in this limit
through a consistent
collection of subtrees, which includes a
tree for every finite set of endpoints in $\R^d$.
Tightness of the distribution, as $\delta \to 0$,
is established for the following two-dimensional examples:
the uniformly random spanning tree on $\delta \Z^2$,
the minimal spanning tree on $\delta \Z^2$
(with random edge lengths), and the Euclidean
minimal spanning tree on a Poisson process of points in $\R^2$
with density $\delta^{-2}$.
In each case, sample trees are proven to have the 
following properties, with probability one with respect to
any of  the limiting measures:
i) there is a single route to infinity (as was known for
$\delta > 0$),
ii) the tree branches are given by curves which are regular
in the sense of H\"older continuity,
iii) the branches are also rough, in the
sense that their Hausdorff dimension exceeds one,
iv) there is a  random  dense subset of $\R^2$, of
dimension strictly between one and two, on the complement of which
(and only there)  the spanning subtrees are unique with
continuous dependence on the endpoints, v)
branching occurs at countably many points in $\R^2$,
and vi) the branching numbers are uniformly bounded.
The results include tightness for the
loop erased random walk (LERW) in two dimensions.
The proofs proceed through the derivation of
scale-invariant power bounds on
the probabilities of repeated crossings of annuli.
\end{abstract}

\vskip .25truecm
\newpage

  \vskip .25truecm         
    \begin{minipage}[t]{\textwidth}
    \tableofcontents
    \end{minipage}
   \vskip .5truecm
\newpage

\masect{Introduction}
\label{sect:intro}

For various systems of many
degrees of freedom, extra insight may be derived by combining
methods of discrete mathematics with considerations
inspired by the continuum limit picture (see e.g.,
\cite{LPS,Car92,Aiz-IMA,Propp,Kenyon}).  
The relation between the
continuum and the discrete perspectives is  through the
{\em scaling limit}.  In this limit the scale on which the
system's defining {\it microscopic\/} variables
can be distinguished is sent to zero,
while focus is kept on features manifested on
a {\em macroscopic scale}.
The first task addressed in this work is a general formulation
of the continuum limit for  stochastic spanning trees.
The existence of limit measures
(which may depend on the choice of subsequence)
is then established for three examples of spanning trees,
all in two dimensions.  The arguments makes use of
the general criteria developed for random systems of
curves  in Ref.~\cite{Aiz-Burch}.
We also derive some basic sample properties of the spanning
trees in the scaling limit.

\masubsect{Three spanning tree processes}

Following are the three examples of random spanning trees
on which we focus in this work.
In each case, the tree connects a set
of sites in $\R^2$ with typical nearest neighbor
distance $\delta \ll 1$.
 \begin{itemize}
\item[UST]
  ({\em Uniformly Random Spanning Tree}) \\
  The vertices to be connected are the sites of the
  regular lattice $\delta \Z^2$, and the spanning tree is drawn
  {\em uniformly at random} from the set of spanning trees
   whose edges connect nearest neighbors in the lattice.
\item[MST]
  ({\em Minimal Spanning Tree}) \\
  The graph is again the
  regular lattice $\delta \Z^2$, with edges connecting nearest
  neighbors. The {\em lengths} associated with the edges
  are determined by {\em call numbers}, which are independent
  identically distributed continuous random variables.
  The spanning tree is the one that minimizes the
    {\em total edge length} (i.e.\ the sum of the call numbers).
\item[EST]
   ({\em Euclidean (Minimal) Spanning Tree})\\
    The vertices of the graph are given by
   a homogeneous Poisson process with density
    $\delta^{-2}$.  We let every pair of vertices be connected
   by an edge whose length is the usual Euclidean
   distance.  The spanning tree is the one that minimizes the total
   edge length.  It may be noted that this spanning tree forms  a subgraph
   of the Voronoi graph of the Poisson process. (In the Voronoi graph,
   a pair of vertices is linked by an edge if and only if there is
   a point in $\R^2$ whose two closest vertices form the given pair.)
\end{itemize}

It is unclear whether our analysis can be extended
to a fourth model, the uniformly random spanning tree on the
Voronoi graph of a Poisson point process.
Such an extension would require a better understanding of
 random walks on the Poisson-Voronoi graph (see the remark
at the end of this introduction).

The scaling limit $\delta\to 0$, can be taken either in fixed
finite regions, $\Lambda \subset \R^d$,  or in conjunction with
the infinite volume limit $\Lambda\nearrow \R^d$.
The analysis of the volume dependence is made easier
by the monotonicity structure
which is discussed here in Section~\ref{sect:FW}.
It is known that for fixed $\delta > 0 $ the limit
$\Lambda \nearrow \R^d$
exists for the spanning trees considered here with either the
free (F) or the wired (W) boundary conditions.  Furthermore,
in any finite dimension the limits coincide for these two
boundary conditions
refs.~\cite{Pem91,Hagg95,BLPS98,Alex-MST,CCN,Alex-RSW}.
The limiting graph,   $\Gamma_{\delta}(\omega)$ 
(with $\omega$ representing the randomness inherent in the
model),   will be free of cycles but in general it need not be
connected
and may instead turn out to be a {\em forest} of infinite trees.

In our analysis of the spanning trees we use the fact that
they can be drawn  with the help of rather efficient algorithms,
employing two processes of independent interest.
The paths of UST obey the statistics of the
{\em loop-erased random walk} (LERW)
\cite{Pem91,Wilson}, while those of MST are related to
the {\em invasion-percolation process\/}~\cite{CCN}.
Through the former correspondence our results return information on the
scaling limit(s) (along subsequences)  of the two dimensional LERW,
which has the same distribution as the path from a predetermined
origin to infinity along the spanning tree (UST).

The relations mentioned above were already employed to shed light on
 the question of unicity of the spanning tree.
Through the relation with the LERW it was shown that
for UST the infinite-volume limit 
a.s.\  consists of a single tree if $d\leq 4$
but of infinitely many trees if $d > 4 $, and  that in any
dimension a.s.\ each tree has a single topological
end (i.e., a single route to infinity) \cite{Pem91,Hagg95,BLPS98}.
As Benjamini and Schramm have observed (private communication)
the situation
in $d=4$ is noteworthy in that in the scaling limit ($\delta =0$)
there will typically be infinitely many trees, while there is
only one tree as long as $\delta >0$.

Less is proven about MST and EST in general, but
it is known \cite{Alex-Mol,Alex-MST,Alex-RSW} (see also \cite{CCN})
that in $d=2$ dimensions $\Gamma_{\delta}(\omega)$
(at $\delta > 0$)
a.s.\ consists of a single tree with a single topological end.
Regarding the upper critical dimension, the situation is less clear.
We think it is possible that the dimension
at  which  the spanning tree is replaced by a forest is $d_c=8$
for MST and EST with non-zero
short-distance cutoff, $\delta > 0$, while the dimension at
which the change occurs for scaling limits
of these models  (i.e., $\delta = 0$) is $d_c=6$.
The heuristics behind the first statement
are discussed in refs.~\cite{New-Stein1,New-Stein2}
in a  context relevant for MST, and essentially the same
heuristics should apply to EST.  The conjecture concerning
the scaling limit is based on the analysis  of percolation
clusters above the upper critical dimension,
discussed in \cite{AizISC}.

\masubsect{Statement of the main results}

Let $\Gamma_\delta (\omega)$ be the infinite-volume limit
of either one of the three spanning tree processes
(UST, MST, or EST)
in $\R^d$, with the ``short-distance cutoff'' $\delta$.
It is an interesting question how to describe the
spanning tree/forest in terms which  remain meaningful in the scaling
limit where the set of vertices becomes dense in $\R^d$.
The approach we take is to describe it through
the collection, denoted below by  $\F_{\delta}(\omega)$,
of all the subtrees spanning finite sets
of vertices.
The benefits are:
\begin{itemize}
\item[i.]
the terminology  makes sense even in
the limit $\delta = 0$;
\item[ii.]
by focusing on the connecting curves and finite subtrees
one can see the tree's ``fractal structure'', which emerges
in its clearest form in the scaling limit;
\item[iii.] the approach can, in principle, be applied in
any dimension.
\end{itemize}

In two dimensions
one could alternatively represent the spanning tree through
its {\em outer contour}, i.e.\ the
line separating it from the dual tree.
The formulation of the scaling limit in terms of such a random
``Peano curve'' was recently suggested by Benjamini
{\em et al.}~\cite{BLPS98}.
Outer contours also play a fundamental role in the broader class
of random cluster models, which includes UST as a limiting
case ($Q\to 0$).  The analysis of such contours played an
important role in physicists' derivation of the exact values
for critical exponents~\cite{Nienhuis,DS}.
(Though not yet rigorously proven, such predictions
appear to be correct. Recent extensions
and applications are discussed in ~\cite{DAA}.)
Let us add, therefore, that our analysis  implies constructive
results also for scaling limits of the outer contours of
the spanning trees studied here.

Thus, we describe a spanning tree/forest by means of the closed
collection of all the subtrees connecting finite collections of
sites. In discussing the infinite volume
limit it is convenient to formulate the curves and trees in
the one-point compactification $\sp^d$ of $\R^d$, which we identify
(via the stereographic projection) with the $d$-dimensional unit sphere.
Since this may result in the blurring of  
the distinction between a spanning tree and a spanning forest, we shall
formulate the difference in Definition~\ref{df:inclusive} below.
Our terminology is built up in the following way
(a more complete discussion of the terms is given in Section 2).

\begin{itemize}
\item[1.]
A  {\em curve  in \/ $\sp^d$} is, for us,  an equivalence
class of continuous functions from the unit interval into $\sp^d$,
 modulo monotone reparametrizations.
Extending this is:

\item[2.] A {\em tree immersed in \/ $\sp^d$}
is an equivalence class of continuous functions from any of the
standard reference trees (see Section 2),
into $\sp^d$.
It will be represented by
the symbol $T^{(N)}(x_1,\ldots, x_N)$, where 
$x_1,\ldots, x_N \in \sp^d$ are the endpoints of the tree.
A subscript $\delta$
may be added to indicate that the tree corresponds to a
model with a short distance cutoff, and a parameter $\omega$
may be added to indicate the random nature of the object.

\end{itemize}

\noindent{\bf Remark:\/}  To avoid confusion let us alert the reader
that for lack of terms, and our reluctance to coin non-intuitive
ones, our terminology may brush against established usage.
Thus, the continuous function defining
an immersed tree need not be invertible, and the intersections which 
occur  need not be transversal,
i.e., the function need not be an immersion in the standard sense.
This notion is natural for our discussion of the scaling limit,
since the trees may have branches
which only appear to intersect, when viewed on the scale of the
continuum, without there being an intersection on the fine scale.

\begin{figure}[htb]
    \begin{center}
    \leavevmode
        \epsfysize=2in
  \epsfbox{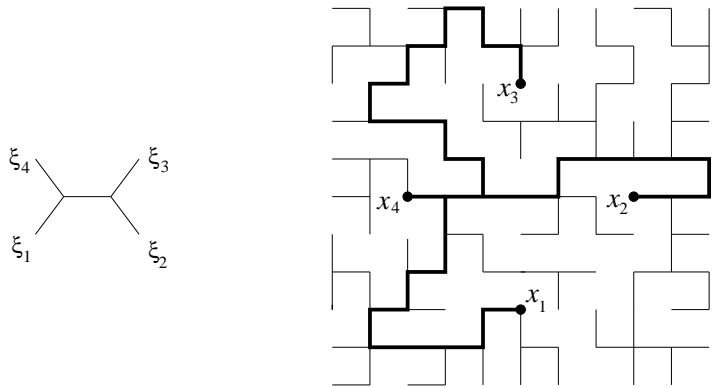}
\caption{\footnotesize
A spanning tree on
a $10\times 10$  grid with free boundary conditions.  Highlighted
is the subtree $T^{(4)}(x_1,\ldots,x_4)$.
The diagram on the left shows a reference tree $\tau$
that can be used to parametrize $T^{(4)}$ (see 
Section~\ref{sect:graphs}.b).
}
\label{fig:subtree}
\end{center} \end{figure} 
\begin{itemize}

\item[3.]  The space of all trees 
immersed in   $\sp^d$ with $N$ endpoints
is denoted here by $\S^{(N)}$.  Note that the restriction of a tree
in $\S^{(N)}$ to  $\R^d$ may be a forest, 
if its branches  pass through infinity.
The spaces $ \S^{(N)}$  are introduced explicitly
in Section~\ref{sect:graphs} along with a metric in
which the distance between two immersed trees
reflects their structure as objects based on curves.
The distance between curves is defined there so that
two curves (or trees) are close if they shadow each other in
a metric on $\R^d$ which  shrinks at infinity. 
Thus convergence in $\S^{(N)}$ means in essence convergence
within bounded subsets of $\R^d$.

\item[4.]  The symbol $\F^{(N)}$ will denote a collection of
immersed trees with $N$ external vertices which forms a
 closed subset of $\S^{(N)}$.  The space of
all such closed collections is $\Omega^{(N)}$.
(Under the induced Hausdorff metric it forms a
complete and separable metric space. )

\end{itemize}

Finally, we are ready to present our full description
of a spanning tree or forest  as a closed collection
of finite trees graded by $N$.

\begin{df}
\begin{itemize}
\item[1.]
A  \underline{spanning forest} for a graph $G$, with
vertices in $\sp^d$ ($d$ fixed at a value which should
be clear from the context),
is represented by a graded collection
$\F = \{ \F^{(N)}  \}_{N\ge 1}$
where: \begin{itemize}
\item[i.] for each $N< \infty$, the collection $\F^{(N)} $
includes a spanning tree
$T^{(N)}(x_1,\ldots, x_n) \in \S^{(N)}$
for each $N$-tuple of vertices of $ G $ ;
\item[ii.] the collection is \underline{inclusive} in the sense that
for any tree $\T \in \F^{(N)}$ (with some $1\le N<\infty$),
all the subtrees of $\T$ are also found in the suitable elements
of the collection;
\item[iii.]  for any two trees, $T_1 \in \F^{(N_1)} $ and
$T_2 \in \F^{(N_2)}$, there is a tree in $\F $
which contains (in the natural sense) both $T_1$ and $T_2$
and has no external vertices beyond those appearing in the
two subtrees.
\end{itemize}

The symbol we use for the
space of all such  collections is $\Om$.
[It forms a closed subset of the product space
${\sf X}_{N\ge 1} \Om^{(N)}$
which we take here with the product topology.]

\item[2.] A spanning forest $\F$ is said to consist of  a
\underline{single spanning tree} in $\R^d$  if every path
$T^{(2)}(x,y) \in \F^{(2)}$ with finite end-points
$x,y \in \R^d$ stays within some finite region of $\R^d$.
[Equivalently (by ii): for every $2 \le  N < \infty$,
each of the immersed trees
$T^{(N)}(x_1,\ldots, x_n) \in \F^{(N)}$
with finite external sites
$\{x_1,\dots , x_N\}\subset \R^d$  is contained in some finite
region of $\R^d$. )

\item[3.] The spanning forest $\F$ is said to be
 \underline{quasilocal} if for any
  bounded region $\Lambda \subset \R^d$ all the trees
 of $\F $ whose external vertices lie in $\Lambda$
 are contained within some
 bounded domain $ \tilde{\Lambda}(\F,\Lambda) \subset \R^d$.
\end{itemize}
\label{df:inclusive}
\end{df}

The probability distribution of
UST, MST, and EST, with the
short distance cutoff $\delta$ as discussed earlier,
correspond to probability measures $\mu_\delta(d \F )$ on $\Om$
(in the appropriate dimension).  
Statements concerning the scaling limits
address limits for the measures $\mu_\delta(d \F )$,
for $\delta=\delta_n \to 0$.
Needless to say, the existence of scaling limits
even along suitable subsequences
is {\it a priori\/}
not obvious since the spaces discussed here
are not even locally compact.  E.g., the tree branches
may, in the limit,  cease to be describable by curves.
Furthermore, in the continuum limit even the most elementary
features could be lost, or appear to be lost:
distinct branches  may  fuse, giving  the appearance of loops
(from the continuum perspective), a tree may turn
into a forest, and multiple paths may open to infinity
(via the stretching to infinity of some of the
connecting paths). In general, concepts which are obvious
or proven for finite graphs need to be re-examined.

Our main results may naturally  be  grouped in two parts.
Following is the first.

\begin{thm}  In $d=2$ dimensions, the following is valid
for each of the
spanning tree processes presented above (UST  and  MST
on $\delta\Z^2$, and EST of density $\delta^{-2}$ on $\R^2$):
\begin{itemize}
\item[i.] {\em (Existence of limit points)}
The collection of measures $\mu_{\delta}(d \F)$
with $0< \delta < 1$ is {\em tight}; every
sequence of $\delta$'s  tending to $0$
includes a subsequence $\delta_n\to  0$ along
which  the  measures
$\mu_{\delta_n}(d \F)$ converge,
in the sense of weak convergence for measures
on the  product space 
${\sf X}_{N\ge 1} \Om^{(N)}$
 to a limit $\mu(d\F)$.
\end{itemize}
For  any of the limiting measures, $\mu$-almost every
spanning forest $\F(\omega)$ 
has the following properties:
\begin{itemize}

\item[ii.] {\em (Locality and basic structure)}
$\F(\omega)$ is quasi-local and
describes a single spanning tree on $\R^2$.

\item[iii.] {\em (Regularity)}
The branches of all the trees in $\F(\omega)$ are
random curves ${\mathcal C}$  with Hausdorff dimensions
bounded above,
\begin{equation}
\dim_{\mathcal H}{\mathcal C} \ \le \ d_{\max}\ ,
\end{equation}
where $d_{\max}<2$ is non-random. 
Furthermore, for any $\alpha<1/2$  
all the curves in $\R^2$ can be simultaneously
parametrized by functions ($g(t)$, $0\le t \le 1$)
which are H\"older continuous of order $\alpha$, i.e.,
each satisfying
\begin{equation}
|g(t) - g(t')| \ \le \ \kappa_{\alpha}(\omega) 
   \left(1+ |g(t)|^2 + |g(t')|^2 \right)
\ |t-t'|^{\alpha} \qquad \mbox{ for all $0\le t < t' \le 1$}
\end{equation}
with  the continuity modulus \
$\kappa_{\alpha}(\omega)$ common to all the branches
of trees in $\F(\omega)$.

\item[iv.] {\em (Roughness)}
Almost surely, all the curves 
(${\mathcal C} \in \F^{(2)}(\omega)$)
are non-rectifiable,
and satisfy also the opposite bound:
\begin{equation}
\dim_{\mathcal H}{\mathcal C} \ \ge \ d_{\min}
\label{eq:dmin}
\end{equation}
with a  non-random $d_{\min} > 1$.
In particular, no branch can be
pa\-ra\-metrized H\"older  continuously with an
exponent  less than $(d_{\min})^{-1}$.
\end{itemize}
\label{thm:main1}
\end{thm}

The convergence asserted  for the measures $\mu_{\delta_n}$
means that
  \begin{equation}
   \int  \psi ( \F) \mu_{\delta_n}( d \F )
   \ \too{n\to \infty} \
    \int  \psi ( \F) \mu( d\F)
  \label{eq:weakconvergence}
  \end{equation}
 for all bounded continuous functions $\psi$ 
 which depend on $\F$ only
  through $ \F^{(N)} $  for some  $N<\infty$
  (for inclusive collections,
  the above is equivalent to permitting dependence on all
   $\{\F^{(1)},\ldots, \F^{(N)} \}$).
This statement may also
be expressed by saying that there is a {\em coupling}, that
 is   a sequence of probability measures $\rho_n$  on
 $\Om\times\Om$  whose marginal
 distributions satisfy
 \begin{equation}
 \rho_n(d \F, \Om) \ = \
  \mu_{\delta_n}(d \F)
 \; ,  \qquad \;
  \rho_n(\Om, d \F) \ = \  \mu(d \F )\ ,
  \label{eq:coupling}
 \end{equation}
  with
  \begin{equation}
 \int_{ \Om \times  \Om}
  \min\left\{1, \dist(\F^{(N)},\F^{\prime (N)}) \right\}
  \; \;  \rho_n(d \F, d {\F^\prime})  \
  \too{n\to \infty} \ 0 \; ,
  \label{eq:vasershtein}
  \end{equation}
  where  $\dist(\cdot, \cdot ) $
  is the Hausdorff distance between closed subsets of
  $ \S^{(N)}$
  based on the metric defined on this space of trees in
  Section~\ref{sect:graphs}.

The proof of Theorem~\ref{thm:main1} utilizes the theory
developed for
systems of random curves in  ref.~\cite{Aiz-Burch}.
The bulk of the analysis consists of the derivation of
the required  criteria, which need to be verified
by model specific arguments.  The criteria amount to scale
invariant bounds on crossing probabilities, which
are presented in the Section~\ref{sect:criteria}.

We believe that the limiting measure $\mu$
of Theorem~\ref{thm:main1} does not depend on the choice
of the subsequence $\delta_n$,  so that
for each of the three processes there is a unique scaling limit.
We further suspect that MST and EST  share a common
scaling limit, based on the
accumulated evidence that the associated critical
percolation models are indistinguishable in this limit,
but that the limit for UST is different.
UST can be presented as corresponding to
the critical Fortuin-Kasteleyn random cluster
model (related to the $Q$-state Potts spin models)
with  $Q \to 0$ along the critical
line~\cite{FK,Dup87,Hagg95}, while MST  is related to critical
percolation, corresponding to $Q = 1$.
The predicted values of characteristic exponents
change with $Q$ (\cite{Nienhuis,DS,Dup87}),
although it should be said that
the exact relation of the exponents of MST with
percolation is not completely clear (to us).

The second set of results describes topological
properties of the spanning trees which emerge
in the scaling limit. To state the results we
need some further terminology.

\begin{df}  For a graded collection of trees $\F \in \Om$
which describes a single spanning tree in $\R^d$:

\begin{itemize}

\item[1.] A point $x\in \R^d$ is said to be a
\underline{point of uniqueness},
if $F^{(2)}$ does not include
a non-constant curve which starts and ends at $x$.

\item[2.] The tree is said to have a
\underline{single route to infinity}
if for any $r>0$ there
is $R(r,\F) <\infty$  such that $\F^{(2)}$ does not
contain a curve spanned by two vertices outside the ball
$B(0;R(r,\F))$ which passes through $B(0;r)$
[i.e., $\infty$ is a point of  uniqueness for $\F$].

\item[3.] $\F$ \underline{branches} at $x\in\R^d$
(and $x$ is called a
\underline{branching point} of $\F$)
if $\F$ includes a tree element
for which $x$ is a vertex of degree at least three,
and the branches meeting at $x$ are non-degenerate in the sense that
they do not collapse to points (i.e., the curves are non-constant).
 
\item[4.]  $\F$ exhibits \underline{pinching}
 at $x\in \R^d$  if  $\F^{(2)}$ includes a curve
which passes through $x$ twice without terminating there.

\end{itemize}
\label{df:point}
\end{df}

It is easy to show (Lemma~\ref{lem:unique})
that if $\F$ represents a single spanning tree in $\R^d$
and $x_1,\dots ,x_N$ are distinct points
of uniqueness, then $\F$ includes
exactly one subtree with external vertices $\eta=\{x_1,\dots ,x_N\}$,
and the corresponding $T^{(N)}$  (viewed as a
tree-valued function of $N$-tuples in $\R^d$)
is continuous at $\eta$.

We prove the following in the scaling limit.

\begin{thm} {\em (Properties of the scaling limits)} \
Let $\mu(d \F)$ be a scaling limit of the measures
$\mu_\delta$ (on $\Om$) discussed in
Theorem~\ref{thm:main1}.  Then $\mu$-almost surely:
\begin{itemize}
\item[i.]
The spanning tree $\F(\omega)$ has a
single route to infinity;
\item[ii.] almost every $x\in \R^2$, in the sense
of Lebesgue measure, is a point of uniqueness for $\F(\omega)$;
\item[iii.]  the set of exceptional points, of non-uniqueness
for $\F(\omega)$, is dense in $\R^2$, and its dimension satisfies
\begin{equation}
2\ >  \ \dim_{\mathcal H} \{ x \in \R^2 |
\ x \mbox{ is not a point of uniqueness for } \F(\omega)
\} \ > \ 1 \; ;
\end{equation}
\item[iv.]  there exists a (non-random) integer $k_o$
so that all non-degenerate trees in $\F(\omega)$
(in the sense that no branches are collapsed to points)
have only vertices of degree less than $k_o$
(see Definition~\ref{df:degree});
\item[v.] the collection of branching points is countable.
\end{itemize}
\label{thm:main2}
\end{thm}

The above assertions follow directly from the power bounds
whose derivation is the main technical part of this paper
(and on which also Theorem~\ref{thm:main1} rests).
In the proof of  Theorem~\ref{thm:main2}
we discuss also a related notion of
the {\em degree} and {\em degree type} of $\F$ at a point $x\in\R^d$ 
(Definition~\ref{df:degree}).

Let us mention that related results were recently presented
for UST by  I. Benjamini~\cite{Ben}, in a work focused on
the large scale features of that spanning tree, seen by ``looking up''
from the lattice scale (while here we focus on the view seen
``looking down'' from the continuum scale).
While the two works, which were carried out independently, differ
in perspectives,  there are
similarities between some of the questions considered and
in the means  employed for their study within the context of UST.

\bigskip\noindent{\bf Remarks}   \ 1) In two dimensions
each spanning tree process has a dual which is also a spanning tree.
Our results for one process imply similar results for the dual, even
without the manifest self-duality which is present in
the case of MST and UST.

2) We expect it also to be true that
in typical configurations of scaling limits of UST, MST and
EST in two dimensions there are no points of branching
of order greater than three, and no points of pinching.
One may approach the proof of such statements
through suitable bounds on the characteristic exponents
(see Section~\ref{sect:conclusion}), however the
analysis presented here does not settle this issue.
A different approach is being suggested by
O. Schramm~\cite{Schramm98}, and a partial result in this
direction (for UST with a short distance ``cutoff'')
can also be found in ref.~\cite{Ben}.

3) An essential ingredient in the  analysis of MST and EST 
is the fact that with positive probability a given point is encircled
on any given scale  by
a critical percolation cluster (see the discussion
after Lemma~\ref{lem:cluster}).
For UST, the corresponding fact is that Brownian
motion in the plane creates loops on all scales
(see the proof of Lemma~\ref{lem:choking}).
The extension of our analysis to the fourth model
mentioned earlier would be facilitated by establishing
that random walks on the Poisson-Voronoi graph resemble
Brownian  motion in that respect,
as stated in the following conjecture (C).
(Some further attention is needed for dealing
with the two sources of randomness: random spanning
trees, in a random graph.)

\noindent{\bf Conjecture (C)}  \ {\em
Let $G(\omega)$ be the random Poisson-Voronoi graph
of density one in $d=2$ dimensions.
For each $x\in \R^2$ and $s\in (0,1)$, let $b_{x,s}(t)$ be the
simple random walk process on   $G(\omega)$ which starts
at the vertex closest to $x$ and continues until
the first exit from the annulus
\[
D_{x,s}= \{ y\in \R^2 \ : \  s|x| \le |y| \le s^{-1}|x| \}\ .
\]
Then there are some $q(s), r_o(s) >0$ such that
 for all starting points with  $|x| \ge r_o(s)$:
\begin{equation}
\P\left( \begin{array}{c}
\text{the trajectory of  \ $b_{x,s}(t)$ separates the} \\
\text{inner and outer boundaries of $D_{x,s}$}
  \end{array} \right) \
\ge q(s) \ > \ 0 \quad .
\end{equation}
(The probability refers here to the double
average corresponding to a random
walk on a random graph.)}

\bigskip

\masubsect{Outline of the paper}

The organization of the work is as follows.
In Section~\ref{sect:graphs} we introduce the space of immersed
trees. Section~\ref{sect:criteria}
contains a summary of the pertinent results from ref.~\cite{Aiz-Burch}.
We recall there two criteria for systems of random curves which permit
to deduce regularity and roughness statements,  as those seen in
Theorem~\ref{thm:main1}.  The criteria require certain
scale-invariant bounds on the probabilities of
multiple traversals of annuli,  and of lengthwise
traversals of rectangles, by curves in  the given random family.
The criteria admit a conformally invariant formulation.
The next two sections present some 
auxiliary results: Section ~\ref{sect:FW} is dedicated to
the very useful free-wired bracketing principle, and
Section \ref{sect:crossing} to preliminary results
on the crossing probabilities for  annuli with various boundary
conditions.  In Section~\ref{sect:reg} we verify
the regularity criterion, treating the three models separately;
in each of the three cases the proof makes use of
a convenient algorithm for generating the tree.
The roughness criterion is verified in Section~\ref{sect:rough}
by means of an argument which applies to all the models discussed here.
In Section~\ref{sect:conclusion}, the results of the previous sections
are combined for the proof of Theorems~\ref{thm:main1}
and~\ref{thm:main2}, followed by
some further comments on the geometry of scaling limits.
The discussion of crossing exponents is supplemented in
the Appendix by deriving a quadratic lower bound
($\lambda(k) \ge \c\ (k-1)^2$) for the
rate of growth of the exponent associated with the
probability of $k$-fold traversals.

\bigskip
\bigskip

\masect{Collections of immersed graphs}
\label{sect:graphs}

Following is the construction of the spaces $\S^{(N)}$
on which we base the description of spanning forests in $\R^d$.
As is mentioned at the end of the section, the concepts
discussed here may be extended to more general
immersed graphs.

\masubsect{Compactification of $\R^d$.}

A convenient way  to encompass in our discussion the infinite
volume limit is to formulate our concepts with the Euclidean metric
replaced by the distance function
$d(u,v)$ defined on $\R^d \times  \R^d$
by
\begin{equation}
d(u,v)\ =\ \inf_\gamma \int_{\gamma}  ds/ (1+|x|^2)\ ,
\label{eq:def-metric}
\end{equation}
where the infimum is over all
continuous paths $\gamma = x(\cdot)$ joining $u$ with $v$, and
$ds$ denotes integration with respect to arclength.
The useful features of the
 metric are:  i) in bounded regions it
is equivalent  to the Euclidean metric, ii) with respect to
it, $\R^d$ is precompact. 
Adding a point at infinity yields the compact space $\sp^d$
which is (via stereographic projection)
isometric to the $d$-dimensional unit sphere.

\masubsect{The space of trees}

For each $N<\infty$  the
space of immersed trees with $N$ external vertices, $\S^{(N)}$,
will be constructed as a union of {\em patches}, each parametrized
by a particular reference tree.  This parametrization is used to
define an initial distance within each patch.
Next, the patches are connected, or sewn together, through an
identification of boundary points, which typically correspond to
trees with some degeneracy.  The space  $\S^{(N)}$ is then metrized
through the imposition of the triangle inequality.

The case $N=2$ corresponds to curves, which
can be defined as equivalence classes of
continuous functions
 $f: [0,1]\rightarrow  \sp^d$, modulo (monotone) reparametrizations.
The distance between two curves, $\C_1$ and $\C_2$, is defined by
\begin{equation}
\dist(\C_1,\C_2)\ :=\
\inf_{\phi_1,\phi_2}\ \sup_ {t\in [0,1]}
d({f_1(\phi_1(t)),f_2(\phi_2(t))}) \; ,
\label{eq:curve_dist}
\end{equation}
where $f_1$ and $f_2$ are particular parametrizations of
$\C_1$ and $\C_2$, and the infimum is over the set of all
monotone (increasing or decreasing)
continuous functions from the unit interval onto itself.

For two curves to be close on $\sp^d$
means that the corresponding curves in $\R^d$ shadow each other
except possibly when they are far from the origin.
Although a Cauchy sequence of  
curves in $\sp^d$ may, in general, converge to
a curve connecting two finite points through infinity,
no such  curves occur in the scaling limits 
of the two-dimensional models discussed here.  
(Systems satisfying
the condition ${\bf H1}$ with $\lambda(2)>0$ 
are easily seen to be quasi-local,
uniformly in $\delta$.)  On the other hand, we do encounter
curves which at one end run off to infinity.

To extend this concept to $N>2$, we replace the interval by a
collection of reference trees.
A {\em reference tree} $\tau$ is a tree graph with finitely
many vertices, labeled as external or internal, with the
external vertices having degree one, and  the
internal vertices having degrees not less than three.
The vertices are connected through links which are
realized as linear continua (intervals) of unit length.
We denote by $N(\tau)$ the number of external vertices.
The number of internal vertices cannot exceed $N(\tau)-2$,
and thus there is a finite catalog of
topologically distinct reference  trees for each given $N<\infty$.

A {\em reparametrization}  of a reference tree $\tau $
is a continuous map $\phi: \tau  \to \tau $
which preserves the sets of internal and external vertices
and is monotone (i.e., order preserving, though not necessarily
strictly monotone) on each link.

\begin{df}
For a given reference tree $\tau$,
a \underline{ tree immersed in $\sp^d$\/ }
indexed by $\tau$ is an equivalence class
of continuous maps $f:\tau\to \sp^d$, with two maps $f_1, f_2$
regarded as equivalent if there are two reparametrizations
$\phi_1$, $\phi_2$ of $\tau$ such that
$f_1 \circ \phi_1 = f_2 \circ \phi_2$.
\label{df:tree}
\end{df}

The collection of immersed trees parametrizable by $\tau$ is
denoted by $\S_{\tau}$, and the collection of
all immersed trees with a given number ($N$) of external vertices
is denoted by $\S^{(N)}=\cup_{\tau: N(\tau) = N} \S_{\tau}$.
Let us note that for each $\tau$ there are elements of $\S_{\tau}$
for which one or more branches
have collapsed to a point (i.e., $f(\cdot)$ is constant on
a link). Such degenerate immersed
trees can be
naturally parametrized by a smaller tree $\tau^\prime$,
and we shall identify it, as an element of $\S^{(N)}$, with a point
in the other collection $\S_{\tau\prime }$.
In this fashion, the set $\S^{(N)}$ may be viewed as
covered by a collection of patches, which are sewn together
and form a connected set.

For each reference tree $\tau$ (with at least two vertices),
a metric $\dist_{\tau}(T_1,T_2)$
is given on $\S_{\tau}$,
by a direct extension of \eq{eq:curve_dist}, in which $\phi_i$
($i=1,2$) denote reparametrizations of $\tau$.
In this metric $\S_{\tau}$ is a complete
separable metric space,
since it is a closed subspace, defined by the incidence relations,
of the space of all $(2 N(\tau) -3)$-tuples of continuous curves
(given by the links -- some of which may be degenerate).

The distance thus defined within each patch yields in a natural way a
metric $\dist(T_1,T_2)$ on  $\S^{(N)}$, defined as the infimum of
the lengths of paths connecting the two points through finite
collections of segments each staying within a single patch.
With this definition $\S^{(N)}$ is a complete separable metric
space, and each $S_{\tau}$ is a closed subspace.

The spaces  $\S^{(N)}$ provide the basic building
element for the space of tree configurations.
As explained in Section~\ref{sect:intro},
we denote by $\Om^{(N)}$ the space of all closed subsets
of $\S^{(N)}$, with the Hausdorff metric, and by
$\Om$ the subspace of  the product
${\sf X}_{N\ge 1} \Om^{(N)}$
consisting of all spanning forests in the sense  of 
Definition~\ref{df:inclusive}.
By construction, $\Om$ is a complete separable metric space.
The following is a useful notion.

\begin{df}  Let $\F\in\Om$ be an inclusive  collection
of trees (see Definition~\ref{df:inclusive})
which represents a single spanning tree for a graph in $\R^d$,    
and  let $T_1,\dots ,T_k$ be a collection of trees in $\F$.  
The trees are
said to be \underline{microscopically} \underline{disjoint} 
if there exists a tree $T$ in $\F$, parametrized as 
$f:\tau \to \dot{\R}^d$,  
which is non-degenerate in the sense that no links are collapsed to
points, and a collection of vertex-disjoint  subtrees 
$\tau_1,\dots ,\tau_k$
of the reference tree $\tau$ so that the restriction
of $f$ to each $\tau_i$  is a parametrization of $T_i$.
\label{df:mic-disjoint} 
\end{df}

Note that our choice of the collections $\F_\delta(\omega)$ guarantees
that for $\delta>0$, microscopical disjointness is equivalent
to disjointness.  In general, microscopically disjoint
subtrees are limits of disjoint subtrees.

\masubsect{Systems of immersed  graphs}

Let us note that the concepts discussed above have a natural
extension to systems of immersed graphs which need not be trees.
Such a generalization may, in fact,  be useful for the description
of the configurations of percolation models (in any dimension).

For the more general system of random graphs one should repeat the
construction in the previous subsection, omitting the requirement
that the graphs which provided the reference index sets $\tau$
be connected and free of loops.
The concepts which would be generalized through this modification
include:
\begin{itemize}

\item[i.] $\S^{(N)}$ --- representing, in the modified definition,
the space of 
 graphs immersed in $\dot\R^d$ with $N$ external
vertices;

\item[ii.] $\Om^{(N)}$ --- the space of closed subsets of $\S^{(N)}$.
\end{itemize}
With this modification
$\F = \{\F^{(N)}\}_{N\ge 1} \in {\sf X}_{N\ge 1} \Om^{(N)}$
represents
a collection of {\em immersed graphs},
to which the notions
of {\em inclusive} configuration and
{\em quasilocal} configuration, introduced
in Definition~\ref{df:inclusive}, also apply.

\masect{Criteria for regularity and roughness}
\label{sect:criteria}

Our proof of Theorem~\ref{thm:main1} employs  the regularity and
roughness criteria developed for systems of random curves
in ref.~\cite{Aiz-Burch}.  Following is a summary of the pertinent
results.
We add here also a brief discussion of the
behavior of the criteria under conformal invariance.
The criteria were derived in the context of a system of
random curves in a finite volume, which in the terminology used
in the present work can be presented as follows.

\begin{df}
A \underline{system of random curves  with a short-distance cutoff
$\delta$} is given by a collection,
$\{ \mu_\delta^{(2)}(d\F^{(2)}) \}_{0< \delta \le \delta_{\max}}$,
 of probability measures
on $\Om^{(2)}$ which provide the probability distributions of
random closed sets of polygonal curves.  The parameter $\delta$
indicates the order of magnitude of the polygonal steps --
in a sense which ought to be clear in the given model.
\label{df:system}
\end{df}

\noindent  {\bf Remarks:}   {\it 1) Motivation.\/}
This terminology is of interest mainly
when there is some
consistency in the formulation of the probability
measures for the different values of $\delta$.
In the examples
considered here these  represent scaled down versions of a
common process, i.e., they are  related by dilations.
The term ``cutoff'' anticipates the possibility that
the measures  $\mu_\delta^{(2)}(d\F^{(2)})$
can be viewed as providing an approximate description of a
process which is defined for $\delta=0$, or possibly some
family of such processes whose approximates  are given by
different sequences with  $\delta_n \to 0$.

\noindent{\it 2) Notation.\/}
The random sets of curves will be  denoted by
 $\F_\delta^{(2)}(\omega)$; and when it be deemed unambiguous
the entire system will be represented by  $\F^{(2)}$, or just
$\F$.  The probabilities evaluated with respect to
$\mu_{\delta}(\cdot)$ will also be referred to as
$\P_{\delta}(\cdot)$.

The possibility raised in Remark 1) requires that the family of
measures either converge to
a limit or at least have accumulation points
as $\delta \to 0$.  Thus the
first question  is one of compactness.  A key
issue here is whether the curves satisfy some uniform regularity
estimates.   A useful tool for the derivation of an affirmative
answer is the general result of ref.~\cite{Aiz-Burch} which
permits to deduce   H\"older continuity bounds
(valid simultaneously for all curves of a  typical configuration
$\F_\delta^{(2)}(\omega)$ in a given compact subset
of $\R^d$) from
estimates on the probability  of multiple  traversals
of a spherical shell.  The required estimate
is formulated as a hypothesis
which needs to be verified  by model-specific arguments.

\masubsect{Regularity criterion}

Denoting the shells by
\begin{equation}
D(x; r, R) =
\left\{y \in \R^{d} \ {\big |}\ r \le |y-x| \le R \right\} \; ,
\label{eq:def-shell}
\end{equation}
and $D(r,R) \equiv D(0; r,R) $,  the required property is stated
as follows.

\begin{itemize}
\item[{\bf (H1)}]
A system of random curves is said to satisfy the hypothesis
{\bf H1} if there is a sequence of exponents
\begin{equation}
 \lambda(k) \too{k\to \infty} \infty
\label{eq:suff-upper}
\end{equation}
such that for each $k < \infty$ and each $s>0$ the crossing
probabilities of spherical shells with radii $0<r<R\le 1$
satisfy
\begin{equation}
\P_{\delta}\left( \begin{array}{c}
       \text{$D(x;r,R)$  is traversed by $k$ vertex-disjoint} \\
\text{segments of a curve in $\F_{\delta}^{(2)}(\omega)$}
\end{array} \right) \ \le\ K(k,s)\,
\left(\frac{r}{R}\right)^{\lambda(k) - s}
\label{eq:lambda}
\end{equation}
uniformly in $\delta\le \delta_o(r,s)$, with some constant
$K(k,s)<\infty$.
 \end{itemize}

It may be noted that $\lambda(1) \le d -1$,  
unless the collection of curves is a.s. empty.
The implication of {\bf H1}  is that with probability one  all the
curves of the random configuration within a compact
set $\Lambda\subset\R^d$ are uniformly equicontinuous, with
a bound that is random but whose distribution does not deteriorate
as $\delta\to 0$. To formulate the result precisely, call a family of
random variables $\kappa_\delta$ {\em stochastically
bounded} as $\delta\to 0$ if
\begin{equation}
\lim_{u\to\infty}\ \sup_{0<\delta\le\delta_o}
\P_\delta\Bigl( |\kappa_\delta(\omega)| \ge u\Bigr)\ =\ 0
\end{equation}
for some $\delta_o>0$.  A family of random variables
$\widetilde \kappa_{\delta}$ is
called {\em stochastically bounded
away from zero}, if the family $(\widetilde\kappa_{\delta})^{-1}$ 
is stochastically bounded.

\begin{thm} {\rm (Regularity and scaling limit  \cite{Aiz-Burch}).} \
Let $\F^{(2)}$ be a system of random curves in a compact region
$\Lambda \subset \R^d$, with short-distance cutoff $\delta$,
and let \ $\{\mu_\delta^{(2)} \}$ be the associated family of probability
measures on $\Om^{(2)}$.  If the system satisfies hypotheses  {\bf H1},
then all the curves $\C \in \F_{\delta}^{(2)}(\omega)$
can be parametrized (through an explicit algorithm)
by functions $f: [0,1]\to \Lambda$
such that for each curve, for all \ $0\le t_{1} < t_{2} \le 1$,
and for every $\eps>0$
\begin{equation}
| f(t_{2}) - f (t_{1}) | \
 \le \kappa_{\eps;\delta}(\omega) \
g(\diam(\C))^{1+\eps} \abs{t_2-t_1}^{\frac{1}{d-\lambda(1)+\eps}} \quad ,
\label{eq:holder}
\end{equation}
with a family of random variables $\kappa_{\eps;\delta}(\omega)$
(common to all $\C\in\F_\delta^{(2)}(\omega)$)
which stays stochastically bounded as $\delta \to 0$. The second factor
depends on the curve's diameter through the function
\begin{equation}
g(r)\ =\ r^{-\frac{\lambda(1)}{d-\lambda(1)}}\quad .
\end{equation}

Moreover, there is a sequence $\delta_{n} \to 0$ for which
the scaling limit
\begin{equation}
\lim_{n\to\infty} \mu_{\delta_n}^{(2)}(d\F^{(2)}) \ := \   \mu^{(2)}(d\F^{(2)})
\label{eq:limit}
\end{equation}
exists, in the sense of (weak) convergence of measures
on $\Om^{(2)}$.
The limit is supported on curves with
\begin{equation}
\dim_{\cal H}(\C)\ \le\ d-\lambda(2)\ ,
\end{equation}
whose parametrization (obtained with the algorithm mentioned
above)  satisfies (\ref{eq:holder}) ---
i.e., it is H\" older continuous with any exponent less
than $1/[d-\lambda(1)]$.
\label{thm:reg-curves}
\end{thm}

\noindent {\bf Remark} Although the above theorem 
was formulated for 
compact subsets $\Lambda\subset \R^d$, the proof requires 
only  that $\Lambda$ is a compact metric space  whose 
Minkowski (box) dimension is at most $d$.  (The H\"older 
continuity condition is to be interpreted in terms of 
the corresponding metric.)   In the present 
work we shall apply it to the Riemann sphere.  

\bigskip

Note that for any spanning tree process
\begin{equation}
\lambda(1) \ =\ 0
\end{equation}
since each point is connected to infinity.  However, we will see that
for UST, MST, and EST, the criterion {\bf H1}
is satisfied on $\R^d$, with
\begin{equation}
\lambda(2) \ > \ 0\
\label{eq:lambda2}
\end{equation}
and $\lambda(k)$ growing at least quadratically with $k$.

\masubsect{Roughness criterion}

The criterion to be verified in order
to prove roughness concerns simultaneous traversals of
cylinders.  We refer by this term to the
solid  body, not its boundary; i.e., a {\em cylinder} of length $L$
and width $\ell$ in $\R^d$ is a set congruent to $I\times B$,
where $I$ is an interval of length $L$,
and $B$ a $(d-1)$-dimensional ball of diameter $\ell$.
A collection of sets $\{A_j\}$ is regarded as {\em well-separated}
if the distance of each set $ A_j $ to
the others is at least twice the diameter of $A_j$.
Following is the hypothesis which is relevant for the study of 
the scaling limit.

\begin{itemize}
\item[{\bf (H2$^*$)}] 
A system of random curves is said here to satisfy the hypothesis
{\bf H2$^*$} if
there exist constants $\sigma\ge 1$, $\rho<1$ and $K < \infty$
such that for every finite collection of well-separated 
cylinders, $A_1,\ldots, A_k $, of
widths $\ell_i$ and lengths
$\sigma\ell_i$ ($i=1,\dots , k$) 
\begin{equation}
\limsup_{\delta \to 0} \ \P_\delta\left( \begin{array}{c}
  \text{each $A_j$ is traversed (``lengthwise'') } \\
  \text{by a curve in $\F_{\delta}^{(2)}(\omega)$ }
  \end{array}
\right)\ \le\ K \rho^k \;  .
\label{eq:suff-lower} 
\end{equation}
\end{itemize}
The asterisk on {\bf H2$^*$} marks a minor modification of  
the condition {\bf H2} formulated in Ref.~\cite{Aiz-Burch}, for 
which  the bound on the probability  
is required to hold for all $\delta < \min_i \ell_i$.   
The pertinent result (which incorporates the comment 
made below) is:
 
\begin{thm} {\rm (Roughness, \cite{Aiz-Burch})} \
Assume that a system of random curves $\F_\delta^{(2)}$
satisfies {\bf H2$^*$}.  Then any measure $\mu^{(2)}$ obtained
as a scaling limit $\delta\to 0$ of the measures $\mu_\delta^{(2)}$
on  $\Om^{(2)}$
is supported on configurations containing only curves
with Hausdorff dimension satisfying
\begin{equation}
\dim_{\cal H} \C\ \ \ge \ d_{\min}
\label{eq:rough-curves}
\end{equation}
with some non-random $d_{\min} > 1$, which depends on
the parameters in {\bf H2$^*$}.
\label{thm:rough-curves}
\end{thm}

\noindent{\bf Remark: }  
Roughness in a random system of curves $\F_\delta^{(2)}$  is 
expressed also on intermediate scales, and it 
does not require any assumption on the existence
of scaling limits.   The full condition {\bf H2} permits to conclude 
lower bounds on the tortuosity of the curves which 
are simultaneously valid on all scales.  
Let $M(\C,\ell)$ be the smallest number
of segments in all the subdivisions of the curve  $\C$ into segments 
of diameters $\le \ell$.
The hypothesis {\bf H2} implies the existence of some 
 $d_{\min} > 1$ such that for any fixed $r > 0$, $ s >
d_{\min}$, and compact $\Lambda\subset\R^d$, the random variables
\begin{equation}
\tilde \kappa_{s,r,\Lambda;\delta}(\omega) \ :=
\ \inf_{ \C \in \F_{\delta;\Lambda}^{(2)}(\omega): 
\ \diam({\cal C}) \ge r }
\ell^s\,M(\C,\ell)
\label{eq:tortuous}
\end{equation}
stay stochastically bounded away from zero, as $\delta \to 0$.
In particular, the minimal number of steps of size $\delta$
needed in order to advance distance $L$ exceeds
$\tilde \kappa \ (L/ \delta)^s$.
This complements Theorem~\ref{thm:reg-curves}, since
under the condition  (\ref{eq:holder}), the random variables
\begin{equation}
\ell^{d-\lambda(1)+\eps}\, M(\C, \ell)
\label{eq:holder-M}
\end{equation}
remain stochastically bounded as $\delta\to 0$.
The general result in \cite{Aiz-Burch}
which implies  both roughness statements is a lower bound on
the {\em capacity} of curves in $\F_\Lambda^{(2)}$. 

One may note that the slightly simpler condition {\bf H2$^*$} 
implies that any scaling limit obeys the full {\bf H2}, 
and thus Theorem~\ref{thm:rough-curves} follows from the statement 
derived in ref.~\cite{Aiz-Burch}.  

For the systems considered here we shall establish the 
hypothesis {\bf H2} in Section~\ref{sect:rough}. 

\masubsect{\  {\bf  H1} under conformal maps}

In discussing   infinite systems it is convenient for us to view
$\sp^d$  as covered by two patches:
the ball $B(R)=\{x\in \R^d \ |\  |x| \le R\}$, with some radius
$R>1$, and the set where $|x| \ge  1/R$.
The inversion  $ ( x \rightarrow  x/|x|^2 ) $ maps the
second patch bijectively onto the  compact region $B(R)$.
The metric defined by (\ref{eq:def-metric})  which we use on $\sp^d$
is invariant under this inversion, and so are the
topologies we defined earlier for the spaces of curves, trees,
and their collections.  
It is useful to know that {\bf H1}
is also stable under inversion:

\begin{lem}
If a system of random curves on $\sp^d$
satisfies the hypothesis {\bf H1}, then so does the system obtained
under the inversion, with the exponents
reduced by not more than a factor of~$2$. Furthermore, if
in the original system the probabilities of simultaneous
$k$ crossings of pairs of disjoint annuli are also bounded by the
products of the corresponding power bounds, then
after the inversion {\bf H1} continues to hold with the
original exponents $\lambda(k)$.
\label{lem:H1-sphere}
\end{lem}

\begin{proof} We need to estimate in powers of $(r/R)$
the crossing probability in
the pre-image of the system of curves in an annulus $D(x;r,R)$.
The  pre-image of any spherical shell
is a set bounded by two spheres (which may degenerate
to hyperplanes, if the boundary of the spherical shell
meets the origin).  Let us
denote the distance between the two spheres as $B$, and
their radii as $\tilde{r}_1 \le \tilde{r}_2$.
We need to distinguish now between two cases:
\begin{itemize}
\item[1)]
if the shell does not include the origin then
the pre-image of $D(x; r,R)$ is compact  --- one of the spheres
encloses the other,
\item[2)] otherwise $( r \le |x| \le R )$, neither of the
two spheres
contains the other, and the pre-image of $D(x;r,R)$ is  the
unbounded set formed by the intersection of their exteriors.
\end{itemize}

In case (1), 
the probability of $k$ traversals in the pre-image of $D(x;r,R)$
is smaller than the probability for the
annulus whose inner boundary
is the smaller of the two spheres and whose outer
radius is ${\widetilde{R}} = \tilde{r}_1+ B$.
Since the system
satisfies {\bf H1} this probability is bounded from above by
$K(k,\eps) (\tilde{r}_1/{\widetilde{R}})^{\lambda(k)-\eps}$
for any $\eps>0$ (see eq.~(\ref{eq:lambda})).

The ratio $(\tilde{r}_1/\widetilde{R})$ may be related to
$(r/R)$ using the invariance of the cross-ratio
 $(z_1 - z_2) (z_3 -z_4)/ [(z_1 - z_3) (z_2 -z_4)]$ 
of the four points at which the surface of $D(x;r,R)$
intersects  the line through $O$ and $x$.  We find:
\begin{equation}
{  (2r) (2R) \over (R+r)^2 } \ = \
{  (2 \tilde{r}_1) (2 \tilde{r}_2) \over
      (2 \tilde{r}_1 + B)(2 \tilde{r}_2 - B)  }  \; .
\end{equation}
It follows that
\begin{equation}
{\tilde{r}_1 \over \widetilde{R}} \ =
\ {\tilde{r}_1 \over \tilde{r}_1 + B}
\ \le \ 4 {r \over R} \; .
\end{equation}
Thus, for such a spherical shell, the image of the system of 
curves under  inversion still satisfies
Eq.~(\ref{eq:lambda}) with the original exponents and constants
$\widetilde K(k,\eps)=4^{\lambda(k)}K(k,\eps)$.

In case (2), the invariance of the cross ratio yields:
  \begin{equation}
 	{ (2r) (2R) \over (R+r)^2} \ =
\ { (2\tilde{r}_1) (2\tilde{r}_2) \over (2 \tilde{r}_1 + B)
(2 \tilde{r}_2 + B) }
  \end{equation}
which implies
\begin{equation}
\left(\tilde{r}_1 \over \tilde{r}_1 + B/2   \right)^2 \ \le \ 
\left(\tilde{r}_1 \over \tilde{r}_1 + B/2   \right)
\left(\tilde{r}_2 \over \tilde{r}_2 + B/2   \right) \ \le
 4 {r \over R} \; .
\end{equation}
To bound the crossing probability in $D(x;r,R)$ we may look
at  two disjoint annuli in the pre-image:
one of inner radius $\tilde{r}_1$ and outer  
 radius $\tilde{r}_1 + B/2$, concentric
with the first ball, and
the other of inner radius $\tilde{r}_2$ and outer radius 
$\tilde{r}_2 + B/2$  concentric with
the second ball.  The {\bf H1}- bound on the crossing probability
within just the first annulus yields for the
image system  the upper bound $ K(k,\eps) (4r/R)^{\lambda(k)/2 -\eps}$.
Under the stronger assumption we recover the full power
$\lambda(k)$.
\end{proof}

It may be interesting to note, though we shall not pursue this
point here,  that the above analysis allows us to deduce that under the
stereographic projection of $\sp^d$ onto the $d$-dimensional sphere,
the Hypothesis {\bf H1}
lifts to {\em conformally invariant bounds}
for the probabilities of $k$ crossings between
pairs of $(d-1)$-dimensional spheres.

\masect{Free-wired bracketing} \label{sect:FW}

The {\em free-wired bracketing principle} is a
useful monotonicity property of both uniform and minimal
random spanning trees, which allows
one to relate the spanning tree on a portion of a large or
infinite graph $G$ to the corresponding object defined
in a subset.  One of its implications
is the existence of the infinite-volume limits with free
as well as with wired boundary conditions.  We shall
encounter other uses below.
In this section we shall briefly recall this known principle
and conclude with a new observation,
expressed here as the {\em free-wired factorization property}, 
which will
be used in the study of the crossing exponents.

Let $G$ be a graph with finite coordination number whose set
of vertices is a locally finite subset
$\V \subset \R^d$, and let $\Lambda \subset \R^d$ be a closed
subset with piecewise smooth boundary
(the reference to such sets is natural in our context, but
it should be clear that the main concepts are not restricted
to graphs immersed in $\R^d$).
The subgraph of $G$ with {\em free boundary conditions}
in $\Lambda$, denoted by $G^F_{\Lambda}$, consists
of the vertex set  $\V_\Lambda = \V\cap\Lambda$, with an
edge between two vertices if and only if there is such an edge in $G$.
Each edge in $G^F_{\Lambda}$ is assigned the length  it had in $G$.

The ``subgraph''
of $G$ with {\em wired boundary} conditions $G^W_{\Lambda}$
is defined similarly, except that rather than simply deleting
all the vertices outside of $\V_\Lambda$, they are merged together
into one vertex $\partial \Lambda$,
called the boundary. In the case of UST and MST, any edge that had
existed between a vertex $x\in \V_\Lambda$ and a vertex $y\not\in
\V_\Lambda$ becomes an edge between $x$ and the
boundary $\partial \Lambda$.
For MST, the corresponding edge length is that of $(x,y)$.
(Note that in $G^W_\Lambda$
there may be more than one edge between a vertex $x$ and 
$\partial\Lambda$
so that $G^W_\Lambda$ is really a multigraph.
In the case of MST, all but the shortest of the multiple
edges joining $x$ to $\partial\Lambda$ may be discarded.)
In the case of EST, the length of the (single)
edge joining a vertex $x\in\Lambda$ with $\partial\Lambda$
is set to equal the Euclidean distance from $x$ to the
geometric boundary of $\Lambda$.

Denote the trees generated by a spanning tree process
on $G^F_\Lambda$ and $G^W_\Lambda$
by $\Gamma^F_\Lambda$ and $\Gamma^W_\Lambda$, respectively, with
$\Gamma^W_\Lambda \backslash \{ \partial \Lambda \}$ the graph obtained
by deleting the special boundary vertex and the edges linking to it.
We slightly abuse the notation by  referring to the restriction
of the tree $\Gamma$ to the subgraph
spanned by the vertices in $\Lambda$ as $\Gamma\cap\Lambda$.

The bracketing principle can be stated as:
\begin{equation}
\Gamma^W_\Lambda \backslash \{ \partial \Lambda \} \  \preceq \
\Gamma \cap \Lambda  \ \preceq  \  \Gamma^F_\Lambda  \; ,
\label{eq:bracketing}
\end{equation}
where $A \preceq B$ means that the set of edges of the random graph
$A$ is {\em stochastically dominated \/ } by the set of edges of $B$,
and where it should be noted that  both the
free and the wired boundary conditions on $\Lambda$ {\em decouple\/}
that region  from the rest of the graph.

The stochastic domination can be expressed through the existence of a
coupling between the two tree processes (in a sense analogous to that seen
in \eq{eq:coupling}) in which a.s. all the edges of $A$ are also contained
in $B$.
For MST, the coupling is provided by constructing spanning trees
simultaneously on $G^F_\Lambda$ and $G^W_\Lambda$ using the same
call numbers (Section~\ref{sect:reg}.b), and for EST by using the same
Poisson points (Section~\ref{sect:reg}.c).  For UST a coupling is
known to exist, though a correspondingly simple explicit coupling
remains unknown.

The bracketing principle implies in particular that the restriction
of the tree $\Gamma^W_\Lambda$ to a fixed ``window'' $\Lambda_o $
is monotone {\em increasing} in $\Lambda $, for $\Lambda \supset
 \Lambda_o$, and that the
similar restriction of $\Gamma^F_\Lambda$ is monotone {\em
decreasing}.  Thus one derives the well-known fact
that the infinite volume limit exists for both free and wired boundary
conditions (separately), and that   the limits
\begin{equation}
\Gamma^{F[W]}(\omega) \ = \ \lim_{\Lambda_n \nearrow {\R}^d}
         \Gamma_{\Lambda_n}^{F[W]}(\omega)
\end{equation}
are independent of the sequence of volumes.
The convergence is in the pointwise sense
 for all three  models under consideration here, provided
the models for the different regions $\Lambda_n$ are
coupled in an appropriate way.
For UST the free-wired bracketing principle appeared implicitly
in~\cite{Pem91}, and was stated and derived explicitly
in~\cite{BLPS98}.  For MST and EST
it appears in ~\cite{New-Stein1,RY}.
It is natural, however, to view  it within the context of
similar principles which have
long been known in related areas; including the early example of the 
Dirichlet-Neumann bracketing for the Laplacian 
(viewed as a quadratic form) 
and the more closely related example of the free-wired bracketing
for the $Q$-state Potts models (discussed for $Q\ge 1$ 
in~\cite {FK,ACCN}).

We shall now add to the collection of monotonicity tools
another useful observation.   Consider the effect of subdividing a 
connected region by a surface which splits it into two sets 
$C$ and $D$, for which we then set the boundary conditions
so that the cutting surface acts (in the natural sense) as a
free boundary for $C$ and as a wired boundary for $D$.
In the interior of $C$ the introduction of the free
boundary along the cut only enhances the spanning tree configuration.
Within $D$ the wired boundary along the cut
diminishes the configuration.  It follows that the original
random spanning tree may be monotonically coupled with either of the 
two separate spanning
tree processes.  We say that the system has the
{\em F/W factorization property} if a simultaneous coupling 
of all these processes can be chosen so that  
the two separate trees in $C$ and $D$ are independent.

\begin{lem} {\em (F/W factorization property)} \
On an arbitrary finite graph, or a finite region in case of EST,
each of the spanning trees considered
here -- UST, MST, and EST, has the free-wired factorization property.
I.e., the three tree processes $\Gamma_{C\cup D}\/$, 
$\Gamma_{C}^{ F}$, and $\Gamma_{D}^{ W}\/ $ can be realized on a 
single probability space so that: 
\begin{itemize}
\item[i.] $\Gamma_{C}^{ F}\/ $ and $\Gamma_{D}^{ W}\/ $
are independent 
 spanning trees (with the indicated boundary conditions along the
 separating surface), 
\item[ii.] within the interior of $\/C\/$,
$ \Gamma_{C}^{ F} \/ $ dominates
$\Gamma\/ $, and 
\item[iii.]  within the interior of $\/D\/ $, $ \Gamma_{D}^{ W} \/$ 
is  dominated by
$\Gamma\/ $.
\end{itemize}
\label{lem:FW}
\end{lem}

\begin{proof}
The existence of such a coupling follows by model specific
arguments.  For MST and EST the argument is most direct,
since the spanning tree is determined by the specified call
numbers in the case of MST, or
specified locations of the points in the case of EST, and the
specified boundary conditions.  For those two cases, the F-W
factorization property is a direct
implication of the F-W bracketing principle and the independence
of the distributions of the variables relevant for the regions
$C$ and $D$.

Another argument is needed for UST.  
As a starting point, we take a coupling
between  the restriction of the full tree 
$\Gamma_{C\cup D}$ to $D$, and 
the ``subtree'' $ \Gamma_{D}^{ W} \/$.
Since $\Gamma_{C\cup D}$ dominates $ \Gamma_{D}^{ W} \/$,
the two measures may
be coupled monotonically, so that claim (iii) holds. 
To construct the coupling with the other component,
$\Gamma_C^{ F}$, we note
that the conditional distribution in
$C$ of $\Gamma_{C \cup D}$, conditioned on its restriction to
$D$,  is just the distribution of UST in $C$ with some partially
wired boundary conditions.  (This is not true for MST,
so the argument makes use of the special structure of UST.)
It follows that the conditional distribution of $\Gamma_{C\cup D}$
within $C$ is always dominated by $\Gamma_C^{ F}$.  
It is therefore possible to extend the measure so that 
(i) and (ii) also hold.  
\end{proof}

In the next section we shall see applications of the above
property.

\bigskip  \bigskip 

\newpage 
\masect{Crossing exponents}
\label{sect:crossing}

This section contains some general considerations regarding
the probability of multiple traversals of spherical shells,
and the exponents $\lambda(k)$ that appear in {\bf H1}.

While $\lambda(k)$  relates to the event that there is
a curve with multiple crossings, we find it useful
to extend the considerations to the events of multiple traversals
by disjoint curve segments -- without requiring those to
 be strung along a common curve.
Thus, modifying slightly the definition of  $\lambda(k)$ given
in \eq{eq:lambda},
we let $\lamst(k)$ be the supremum of
all exponents $s$ such that, for all spherical shells with radii
$0<r<R\le 1$,
\begin{equation}
\P_{\delta}\left( \begin{array}{c}
       \text{$D(x;r,R)$  is traversed by $k$}\\
\text{microscopically disjoint curves in $\F_\delta^{(2)}(\omega)$} \\
\end{array} \right) \ \le\ K(k,s)\, \left(\frac{r}{R}\right)^{s}
\label{eq:lambda-bar}
\end{equation}
holds uniformly in $\delta\le \delta_o(r,s)$ with some constant
$K(k,s)$.
Since we relaxed here the condition seen in \eq{eq:lambda},
the exponents are related by
\begin{equation}
\lambda(k)\ \ge\ \lamst(k)\ .
\label{eq:lambdabar}
\end{equation}
The regularity assumption {\bf H1} will be verified by
establishing lower bounds on $\lamst(k)$.

In our discussion we shall make use of the
free-wired bracketing principle and the
F/W factorization property.
The results
of this section hold for any random  spanning tree model
to which these principles applies, regardless of the dimension.

\bigskip

\masubsect{The exponents $\phi(k)$, $\gamma(k)$, and the
geometric-decay property}

In the study of the exponents it convenient to introduce two
additional variants, which correspond to the crossing
probabilities with different combinations of boundary conditions.
The boundary conditions are indicated here in the superscript.
For example,
the graph $G^{F,W}_{r,R}$ is defined by placing on $D(r,R)$
the {\em free} boundary conditions at $r$ and the
{\em wired} boundary conditions at $R$,
i.e., \ deleting the vertices inside $B(r)$ and outside $B(R)$,
and adding a single vertex to the graph representing $\partial B(R)$.

\begin{figure}[htb]
    \begin{center}
    \leavevmode
        \epsfysize=2in
  \epsfbox{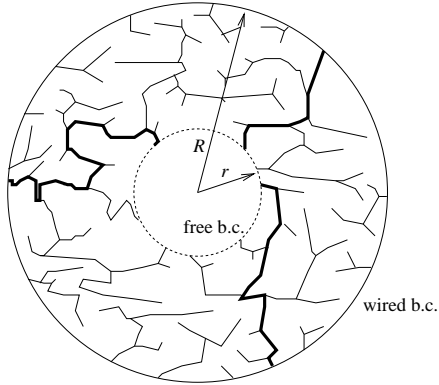}
\caption{\footnotesize
The tree depicted here has $k=3$ disjoint crossings of the
annulus $D(x;r,R)$ with free-wired boundary  conditions.
Note that each point is connected to the wired
boundary by a unique path, while there are many paths to the free
boundary. The exponents $\gamma(k)$ appear in bounds for the
probability of such $k$-crossing events.}
\label{fig:gamma}
\end{center}
\end{figure}

Since traversing means reaching the boundary, or beyond,
some adjustment in the definition is needed at  the
free boundary.  We do that by defining boundary
sites, and then saying that a path along the edges of
$G^{F,W}_{r,R}$  {\em traverses} $D(r,R)$ if
it connects a vertex on the free boundary at $r$ with
a vertex on the wired boundary at $R$.
In the case of the lattice models (UST and MST),
a vertex $x$ in $D(r,R)$ is said to lie {\em on  the free boundary}
of $G^{F,W}_{r,R}$ at $r$ if the original graph $G$ contains an
edge joining $x$ to a vertex inside the ball $B(r)$. In the
case of EST  the defining condition is that
the Voronoi cell of $x$ touches the Voronoi cell
of $\partial B(r)$, or equivalently,
that there exists a disc which intersects $B(r)$ and contains
$x$ but no other vertex of $G$.  The free boundary at $R$ is
defined analogously.
 Two traversals are {\em disjoint}, if they do not share
 any vertices.  
(The alternative definition, based on edge disjointness would 
result in the same exponents.)     

We now define two new families of exponents which play an 
auxiliary role.   Let $\phi(k)$ be the  supremum of all
$s>0$ such that for every shell $D(r,R)$
\begin{equation}
\P_\delta\left(
\begin{array}{c}
        \text{$\Gamma^{F,F}_{r,R}$ includes at least}  \\
        \text{$k$ disjoint traversals of $D(r,R)$}
\end{array}
  \right)  \ \le\ K(k,s) \left( \frac{r}{R} \right)^s
\label{eq:phi}
\end{equation}
for $\delta\le \delta_o(r,s)$, with some constant $K(k,s)$ which
does not depend on $\delta$.

Similarly,  let  $\gamma(k)$ be defined by the condition
\begin{equation}
\P_\delta\left(
\begin{array}{c}
        \text{$\Gamma^{F,W}_{r,R}$ [$\Gamma^{W,F}_{r,R}$] includes at
least}  \\
        \text{$k$  disjoint traversals of $D(r,R)$}
\end{array}
  \right)  \ \le\ K(k,s) \left( \frac{r}{R} \right)^s
\label{eq:gamma}
\end{equation}
interpreted as above (with independently defined constants).
In \eq{eq:gamma}  it is required that the bound holds
for both mixed boundary conditions.

All three families of
crossing exponents are clearly nondecreasing with $k$.
We expect that  $\lambda(k)=\lamst(k)=\phi(k)=\gamma(k)$.
It is shown below that
\begin{equation}
\lamst(k) \ \geq \ \phi(k) \ \geq \
\gamma\left(\left\lceil\frac{k+1}{2}\right\rceil\right)\ ;
\end{equation}
free-wired bracketing easily implies that $\gamma(k)\ge\phi(k)$.

The desired statement: $\lamst(k)\to \infty$, will be derived
by showing that in the UST, MST, and EST models
the crossing probabilities have the following {\em geometric-decay}
property (in $k$) for shells of fixed aspect ratio: There exist
constants $s>0$ and $\sigma>1$ so that
\begin{equation}
\P_\delta \left(\begin{array}{c}
\text{$\Gamma^{F,W}_{r,R} [\Gamma^{W,F}_{r,R}]$ contains at least $k$ }\\
\text{disjoint traversals of $D(r,R)$}\end{array}\right)
\ \le  \left( \frac{r}{R} \right)^{s(k-1)}\
\label{eq:geom-decay}
\end{equation}
holds for all $R \ge \sigma r$, provided $\delta\le\delta_o(r)$.
This  implies:
\begin{equation}
\gamma(k) \ \ge \ s\ (k-1) \; ,
\end{equation}
which suffices for our main purpose.
However, note that \eq{eq:geom-decay} also implies more,
since our definition of the exponents left room for some prefactors,
i.e., it concerned only the asymptotic behavior of the crossing
probability as
$R/r\to\infty$, at fixed $k$.  In the appendix we show, by an argument
of more general applicability which uses the geometric-decay property,
that the actual rate of growth of the exponents is even higher, with
\begin{equation}
\gamma(k) \  \ge  \ \beta(k-1)^2
\label{eq:gamma-quad}
\end{equation}
with some $\beta > 0$.

\masubsect{Comparison of the exponents}

\begin{lem}
$\lamst(k), \gamma(k) \geq \phi(k)$.
\label{lem:lambda-phi}
\end{lem}

\begin{proof}
Recall that $\lamst(k)$ pertains to events
involving a single tree containing
multiple traversals of a spherical shell $D(r,R)$.
Imposing free boundary conditions on
the inner and outer boundaries  of the shell
is a monotone operation which
preserves the traversals.
Thus, the first claim seems to be an immediate consequence of the
bracketing principle (\ref{eq:bracketing}). There is however
one scenario which requires a bit more attention:
Some of the traversals (appearing in the definition (\ref{eq:lambda-bar})
of $\lamst(k)$) may be realized
by an edge which crosses  the annulus $D(r,R)$ without ``stepping'' on
a point in it. In the Poisson-Voronoi graph, the one case
in which this warrants some attention, this
event can occur only if  within the
region $D(r,R)$, there is a disc of diameter at least
$(R-r)$ which contains no Poisson points.  The probability of that is
not greater than approximately $e^{-\c\ (R-r)^2/\delta^2}$.  Such a correction
term plays a negligible role and does not interfere with our ability to
conclude that $\lamst(k) \geq \phi(k)$.

The second claim, $\gamma(k)\ge\phi(k)$,  follows directly
from the free-wired bracketing principle.
\end{proof}


\begin{lem}
$\phi(k) \geq \gamma(\lceil (k+1)/2\rceil)$.
\label{lem:phi-gamma}
\end{lem}

\begin{proof}
Assume $\Gamma^{F,F}_{r,R}$ contains
$k$ (or more) disjoint paths traversing $D(r,R)$.
Label the traversing curves such that $\C_i$
connects a point $p_i$ on the free boundary
at $r$ to a point $q_i$ on the free boundary at $R$. Let $T$ be the subtree
of $\Gamma^F_{r,R}$ spanned by the points $p_1,q_1,\dots ,p_k,q_k$;
it consists of $\C_1, \dots , \C_k$ and $k-1$ ``joining paths''.

Divide $D(r,R)$ into $m$ subshells of aspect ratio $(R/r)^{1/m}$.
Assume that the $p_i$ lie in the innermost,
and the $q_i$ in the outermost subshells --
for UST and MST this happens with certainty if $\delta\le\delta_o(r,m)$,
and for EST the probability of it failing introduces
a negligible correction which is exponentially
small in $\delta^{-2}$, as discussed above.

\begin{figure}[htb]
    \begin{center}
    \leavevmode
        \epsfysize=2in
  \epsfbox{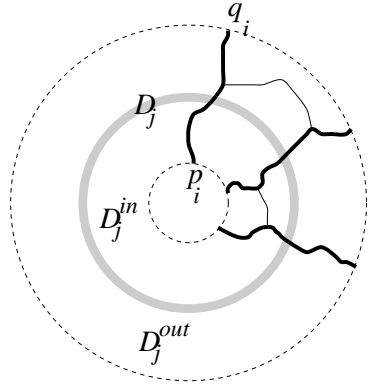}
\caption{\footnotesize
A subtree consisting of $k=3$ disjoint traversals
and $k-1=2$ joining curves.
}
\label{fig:phi-gamma}
\end{center}
\end{figure}

Wiring both boundaries of a subshell $D_j$
divides $D(r,R)$ into an inner shell $D^{\text{in}}_j$
(with free-wired boundary conditions), an outer shell $D^{\text{out}}_j$
(with wired-free boundary conditions),
and the middle (wired) subshell $D_j$.   It is possible to choose $j$
such that each of $D^{\text{in}}_j$ and $D_j^{\text{out}}$ contains at most
$(k-1)/2$
of the joining paths in $T$. With this choice, the intersection
of $T$ with $D^{\text{in}}_j$ consists of at least
$(k+1)/2$ disconnected subtrees, each of which contains at least
one of the points $p_i$, and hence a
traversal of $D^{\text{in}}_j$.  Since wiring the middle shell only
suppresses edges, but
each $p_i$ remains connected to the wired boundary of $D_j$, there are
$(k+1)/2$ traversals of $D^{\text{in}}_j$ with free-wired boundary conditions.
(If $(k+1)/2$ is a half-integer, we may round up.)
By the same reasoning, there are at least $(k+1)/2$  paths
traversing the outer shell $D^{\text{out}}_j$ with wired-free boundary
conditions. Summing over the
possible positions  of $D_j$, and using  the independence
of the tree processes on
$D^{\text{in}}_j$ and $D^{\text{out}}_j$, we find that for
each $s<\gamma(\lceil (k+1)/2\rceil)$
(see the definition of $\gamma(k)$ in (\ref{eq:gamma})) we have
\begin{eqnarray}
\P_\delta\left(
\begin{array}{c}
        \text{$\Gamma^{F,F}_{r,R}$ includes $k$ disjoint}  \\
        \text{traversals of $D(r,R)$}
\end{array}
\right)
 &\le& \\
 \nonumber
  &\le&
 m\ \Bigl[K(\lceil (k+1)/2\rceil,s)\Bigr]^2\
          \left( \frac{r}{R} \right)^{s(1-1/m)}
 + E_\delta\ ,
\end{eqnarray}
where $E_\delta$ is a correction term of order
$O(e^{-\c\ (R-r)^2/\delta^2})$. Since $m$ was arbitrary,
it follows that $\phi(k) \geq \gamma(\lceil (k+1)/2\rceil)$.
\end{proof}


\masubsect{A telescopic bound}

A very useful consequence of
the F/W factorization property is a telescopic bound
of the crossing probabilities, which is expressed in the
following lemma.   It yields lower bounds on the exponents
 $\gamma(k)$ from bounds on the crossing probabilities of
 spherical shells with a fixed aspect ratio.

\begin{lem} {\em (Telescopic principle)} \
For each of the spanning trees considered here
(UST, MST, and EST), and in any dimension,
the following is satisfied for any \
$r_1 < r_2 <\cdots  <r_m$ \ and any integer $k$
\begin{eqnarray}
\P_\delta\left(
\begin{array}{c}
        \text{$\Gamma^{F,W}_{r_1,r_m}$ contains $k$ disjoint}  \\
        \text{traversals of $D(r_1,r_m)$}
\end{array}
  \right) &\leq&  \\
&& \nonumber
 \hskip -7cm    \le \
   \prod_{j=1}^{m-1} \biggl[ \P_\delta\left(
\begin{array}{c}
        \text{$\Gamma^{F,W}_{r_j,r_{j+1}}$ contains $k$ disjoint}  \\
        \text{traversals of $D(r_j,r_{j+1})$}
\end{array}  \right)
      \;+\; \P_\delta\left(
\begin{array}{c}
        \text{$D(r_j,r_{j+1})$ is crossed}  \\
        \text{by an edge in $\Gamma^{F,W}_{r_1,r_m}$  }
\end{array} \right)
 \biggr].
\label{eq:telescopic}
\end{eqnarray}
The analogous relations are also valid for the free-free and
the wired-free boundary conditions.
\label{lem:telescopic}
\end{lem}

\remark As mentioned before, the possibility of a ``long edge''
introduces  a correction (the second term on the right)
  whose effect on the exponents
discussed here is negligible.

\begin{proof}
Consider the effect of subdividing a
spherical shell $D(r,R)$ by a
sphere of radius $\tilde{r}$, with the
boundary conditions placed so that the cutting surface
acts as a free boundary for the outer shell $D(\tilde r, R)$,
and as a wired boundary for the inner shell $D(r,\tilde r)$
(so that we end with free-wired boundary conditions
on each subshell).  As we saw in Lemma~\ref{lem:FW} there
 exist  a coupling between
the spanning tree in $D(r,R)$ and the product measure of the
spanning trees in the subshells,
which is separately monotone in the two regions.
On the outer subshell, introducing the free
boundary along the cut  only enhances the configuration.
On the inner subshell, introducing the wired boundary along the cut
diminishes the configuration; however, even in
the diminished spanning tree,
each site remains connected to the
wired boundary.  It follows that every traversal
of $D(r,R)$ of the original configuration  which contains at least
one vertex in each subshell is preserved  as a traversal of both
subshells in the final configuration.
The independence of the two
components, up to the correction which was mentioned
explicitly above, implies the statement for
$m=2$.  The rest is by induction.
\end{proof}

\masubsect{Extension of the bounds to $\delta = 0$}

Another  important property of the exponents, which is valid
in a great deal of generality, is their ``lower semicontinuity'',
in the following sense.

\begin{thm}
Let $\{ \mu_\delta( d\F)  \}$ be a system of random trees
with a short distance cutoff $0 < \delta \le 1$, for which
some of the exponents $\lamst(k)$, $\lambda(k)$ and
$\gamma(k) \geq \phi(k)$  have strictly positive values.
Then the corresponding upper bounds, expressed by equations
(\ref{eq:lambda}), (\ref{eq:lambda-bar}),
(\ref{eq:phi}), and (\ref{eq:gamma}),
continue to apply also at $\delta = 0$
for any limiting measure
$ \mu(\cdot) \ = \ \lim_{\delta_n \to 0} \mu_{\delta_n}(\cdot) $  
(with respect to weak convergence of probability measures on $\Om$).
To be explicit:
the above hold with
unchanged values of the exponents
$\lambda(k)$, \ldots ,
though the optimal exponent values for $\mu$ (at $\delta = 0$)
may be even greater.
\label{thm:semicontinuity}
\end{thm}
\begin{proof}
It is convenient to carry out the argument using the coupling
formulation of convergence, as in \eq{eq:vasershtein}
(with
the distance function evaluated between the finite volume
configurations $\F_{\Lambda}^{(N)}$).
Let us first note that for each given annulus, or
spherical shell, the set of tree configurations which satisfy
the corresponding
multiple crossing condition forms a closed subset of $\Om$.
Therefore its measure under $\mu_{\delta}$ would be
upper semicontinuous, i.e., upward  jumps (as $\delta\to 0$)
are not excluded.
Such discontinuities occur if the
approximating configurations exhibit curves which
stretch and span $D(r,R)$ in the limit.  The probability of
that can be bounded by the crossing events of the arbitrarily
narrower shells (or annuli) $D(r+\eps,R-\eps)$.
This correction can be easily incorporated into the
optimization parameter $s$; the result being that the
upper bounds continue to hold with the $\delta>0$ value of
the exponents $\lamst(k)$, \ldots, $\phi(k)$.
\end{proof}

\masect{Verification of H1 in two dimensions}
\label{sect:reg}

We verify the regularity criterion {\bf H1}
for the three models separately, by reducing it
in each case to a property
of a well-studied random model. Specifically,
for UST, we refer to known properties of random walks,
and for MST and EST to properties
of two independent  percolation processes.
Unlike the previous section, the discussion is
now narrowed to $d=2$.

\masubsect{Uniformly random spanning tree}

We find it useful to construct UST
with the loop-erased random walk algorithm
(\cite{Wilson}).  The {\it current tree\/} starts
out consisting of a single vertex, called
the {\em root}.  The algorithm runs loop-erased
random walk (LERW), starting from any vertex, until
the current tree is
reached.  At that point, the loop-erased
trajectory is added to the
current tree.  This process continues until all vertices have been
adjoined to the tree, which is then uniformly random,
regardless of the choices of the root and the starting
points for the LERW's.

\begin{lem}
Consider UST on an annulus $D=D(r,R=3r)$ with
any (e.g.\ free-free, free-wired, or wired-free) boundary conditions, 
and
let $T$ be a connected subtree (of the appropriate graph for those
boundary conditions)
containing at least one traversal of the annulus.
Condition upon the
edges of $T$ being contained in UST.
Then except with probability $3^{-\alpha}$
(uniformly in $\delta\le\delta_o(r)$, with some $\alpha>0$ which
does not depend on  $r$ or $T$),
UST contains also a {\it choking surface},
which is a collection of vertices that are connected within the 
spanning tree to $T$
via paths that stay within the
annulus (i.e.\ avoid the boundaries), and such that every path
crossing the annulus intersects the choking surface.
\label{lem:choking}
\end{lem}

\begin{proof}
To pick a random spanning tree conditioned to contain some set of
edges (in this case the edges of $T$), we can contract the given
edges, and take the remaining edges from a random spanning tree of the
contracted graph.  Since by assumption $T$ is connected, we can
implicitly contract the edges of $T$ by initializing the current tree
to be $T$ and build up the rest of the tree via loop-erased random
walks.  Let $x$ be a point approximately at radius $2 r$ (i.e., far
from both the inner and outer boundaries of the annulus).  Let
$x_1=x,x_2,x_3,\ldots,x_n$ be the vertices which a random walk (unobstructed
by $T$) visits, up to and including the time that either (1) it hits a
boundary, or (2) its loop erasure makes a non-contractible loop, i.e.\
the loop-erasure of $x_1,\ldots,x_{n-1}$ together with the edge
$(x_{n-1},x_n)$ includes a loop $\hat{C}$ winding around the inner circle.
Recall that we start with $T$ as the current tree.  When we build the
random spanning tree containing $T$, the first $n$ ``choices'' that we
make will be $x_1,x_2,\ldots,x_n$, in that the choices of where to
start the loop-erased trajectories, and the random choices of where
the trajectories go are, are determined by $x_1,x_2,\ldots,x_n$.  I.e., the
first segment adjoined to the current tree is the loop-erasure of
$x_1,x_2,\ldots,x_i$, where $x_i$ is the first vertex from the
sequence already in the current tree $T$.  The second segment adjoined
to the tree is the loop erasure of $x_{i+1},x_{i+2},\ldots,x_j$, where
$x_j$ is the first vertex in the rest of the original RW sequence that
is in the current tree at that point.  We continue in this fashion; if
constructing the tree requires more choices (steps) after the first
$n$, then these are drawn from fresh coin flips.

Consider the random walk winding event described above
(i.e., that (2) occurs before (1)), and let $C$
denote the set of vertices in the noncontractible cycle $\hat {C}$.
We claim that $C$ is contained in the current tree by step $n$, and
comprises a choking surface. To see this, note first that by planarity
$C$ meets every path connecting the inner and outer boundaries.
Secondly note that every loop that is erased in the construction of
the spanning tree by step $n-1$ must necessarily also be erased from
the loop-erasure of $x_1,\ldots,x_{n-1}$.  (This takes a moment's
thought, and the ``cycle-popping'' viewpoint of the LERW construction
(\cite{Wilson}) may help.)  Thus at step $n-1$ each vertex in the
cycle $\hat {C}$ is either contained in the current tree or the
current loop-erased trajectory.  In particular, $x_n$ (visited at a
previous time step) is in the current tree, since the cycle $\hat {C}$
intersects the crossing of the annulus contained in $T$, and the
portion of $\hat {C}$ prior to this intersection will not be in the
current loop-erased trajectory.  When the walk again reaches $x_n$ at
step $n$, all the vertices in the current loop-erased trajectory are
added to the tree.  Since the walk never visited either boundary, each
vertex in $C$ is connected to the initial current tree $T$ via a path
that avoids the boundaries.

It follows from a standard fact about Brownian motion that there is
some positive number $p$ so that whenever $\delta\le\delta_o(r)$, with
probability at least $p$ the loop-erased random walk started from
point $x$, if it is unobstructed by $T$, will wind around the inner
circle and intersect itself before reaching either boundary.  The
assertion follows by choosing $\alpha$ so that $3^{-\alpha}=1-p$.
\end{proof}

\begin{cor} Let $\alpha$ be the exponent of Lemma~\ref{lem:choking}.
UST has the geometric-decay property (\ref{eq:geom-decay}) with
$s=\alpha$ on shells with mixed boundary conditions (free-wired or
wired-free) and aspect ratio $3$.
\label{cor:geom-decay-UST}
\end{cor}

We remark that this corollary is essentially contained in
the proof of part 2 of Theorem 2 of Benjamini's article (\cite{Ben}).

\begin{proof}
We can construct the spanning tree on the spherical shell by starting
LERW's along each point on the free boundary, and only after all the
free boundary vertices are in the current spanning tree, start the
LERWs at other vertices.  Suppose that the LERW from some vertex on
the free boundary makes it to the wired boundary, making the
$k$th ($k\geq 1$) disjoint traversal of the annulus.  We can upper
bound by $3^{-\alpha}$ the
probability that eventually there is a $(k+1)$st disjoint traversal:
By Lemma~\ref{lem:choking}, with probability at least $1-3^{-\alpha}$
there is a choking surface relative to the tree built so far.  But the
tree built so far has only $k$ connections to the wired boundary, so
each vertex on the choking surface is connected to the wired boundary
along one of these $k$ connections.  A
$(k+1)$st traversal disjoint from the previous $k$ traversals would add
a second path from the wired boundary to some vertex in the choking
surface.
\end{proof}

\begin{cor} {\rm({\bf H1} for UST)} For all $k\ge 1$,
\begin{equation}
\gamma(k+1)\  \ge \ \gamma(k) + \alpha \quad ,
\label{eq:lin-gamma-UST}
\end{equation}
with $\alpha>0$ as in Lemma~\ref{lem:choking}. In particular,
{\bf H1} holds for UST with
\begin{equation}
\lambda(k) \ \ge \lamst(k)\ \ge \ \phi(k)\ \ge\ \frac{\alpha}{2}(k-1) \   .
\label{eq:lin-lambda-UST}
\end{equation}
\label{cor:H1-UST}
\end{cor}

\begin{proof}  The first claim is an immediate consequence
of Corollary~\ref{cor:geom-decay-UST} and Lemma~\ref{lem:telescopic}; the
second claim also  uses Lemma~\ref{lem:lambda-phi} and~\ref{lem:phi-gamma}.
\end{proof}

\masubsect{Minimal spanning tree}

The arguments in this subsection are based on the relation between
MST and  critical Bernoulli percolation.
We begin with the natural coupling between the two processes.

Let $\{u_b\}$ (indexed by the edges $b=\{x,y\}$ in $\delta\Z^d$) be
a family of independent
random variables which are uniformly distributed
on $[0,1]$.  These  are the {\em call numbers} which determine the edge
lengths mentioned in the introduction; they
already give a coupling  to Bernoulli  percolation
for all parameter values $p$, i.e., a way to realize the models for
different $p$'s on the same probability space.
To realize Bernoulli percolation for a parameter value $p\in [0,1]$,
we simply call an edge $b$ $p$-occupied if $u_b<p$; then the
$p$-occupied edges (and their associate $p$-clusters, $p$-paths, etc.) are
a realization of density-$p$ Bernoulli percolation.

For given values of the call numbers,
MST can be constructed  on a bounded region $\Lambda$ by
the following {\em invasion} process.
Starting with any vertex as the root, the tree grows by
adding at each step the neighboring edge with the lowest call
number, provided no loop (and no loop through a wired boundary)
is formed; if a loop would be formed, the edge is discarded.
The construction terminates when the tree spans all vertices.
The result is the unique edge-length minimizing spanning tree,
regardless of the choice of the root,
provided that no two edges were assigned the same call number.

\begin{lem} \label{lem:cluster} For MST on a finite graph,
in any dimension and with any of the boundary conditions used here,
if an edge $b$ is vacant in a
configuration, then almost surely its endpoints are connected with each other
(possibly through a wired boundary) by a $p=u_b$-path.
\end{lem}

\begin{proof}
Construct the tree as described above, with one of the endpoints
of $b$ as the root. With probability one, all edges other
than $b$ have call numbers different from $u_b$.  If $b$ is vacant,
then the subtree connects the root to the second endpoint
of $b$ using only edges with call numbers less than~$u_b$.
\end{proof}

Denote by $p_c=p_c(d)$  the Bernoulli percolation
critical value (which for $d=2$  is $p_c=1/2$ \cite{Kesten}).
It is an implication of the Russo-Seymour-Welsh
theory \cite{Russo,SW} that for critical Bernoulli
percolation in $\delta\Z^2$
\begin{equation}
\P_\delta \left( \begin{array}{c}
\text{$D(r,3r)$ is traversed }\\
\text{by a $p_c$-path} \end{array} \right)\ \le\
3^{-\alpha} \quad (0<\delta\le\delta_o(r))\ ,
\label{eq:def-alpha-MST}
\end{equation}
with some $\alpha>0$.
Bounding the  probability of crossings by disjoint
$p_c$-paths by using the van den Berg-Kesten inequality~\cite{BK}
results in the geometric-decay property that
\begin{equation}
\P_\delta \left(\begin{array}{c}
\text{$D(r,3r)$ is traversed by}\\
\text{at least $k$ disjoint $p_c$-paths}\end{array}\right)
\ \le\ 3^{-\alpha k}\ .
\label{eq:geom-decay-perc}
\end{equation}
Since spatially separated events are independent,
a telescopic argument analogous to Lemma~\ref{lem:telescopic}
implies that {\bf H1} holds for critical Bernoulli percolation
in two dimensions, with exponents  $\gamma_{B}(k)$ satisfying
\begin{equation}
\gamma_{B}(k)\ \ge\ \alpha k\, \ > \  0\  ,\quad (0<\delta<\delta_o(r))\ .
\label{eq:H1-MST}
\end{equation}
Note that the probability of crossing
events for Bernoulli percolation does
not depend on the boundary conditions placed on $D(r,R)$.

The object is to bound the probability that
MST  contains $k$ paths traversing an annulus
in terms of related events in critical Bernoulli percolation.
For a given locally finite connected graph $G$ embedded in the
plane, consider the {\em dual} graph $G^*$. Its vertices are the cells of
$G$ (i.e., the connected
components of the complement in $\R^2$ of the union
of the embedded edges of $G$).   There is a dual edge $b^*$ joining two dual
vertices for each common edge in the boundary of the corresponding two cells.
In general, $G^{**}=G$, but note that $G^*$ can be a multigraph.
In particular, $\delta\Z^2$ can be
drawn with vertex set $\delta\Z^{2*}=
\delta\Z^2+\left(\frac{\delta}{2},\frac{\delta}{2}\right)$,
and each dual edge $b^*=\{x^*, y^*\}$ is the perpendicular bisector of
some edge $b=\{x,y\}$.
The dual of the graph $G^{F,W}_{r,R}$ contains
a single vertex $\partial B(r)^*$ dual to the cell
inside the free boundary at $r$ which plays the role
of a wired boundary for $G^{F,W*}_{r,R}$. A row of vertices
dual to the cells touching the wired-in  point $\partial B(R)$
plays the role of a free boundary for the dual.  The
analogous description holds for $G^{W,F*}_{r,R}$, with the roles
of the boundaries at $R$ and $r$ interchanged.

A dual bond $b^*$ is called $p$-occupied when $b$ is $p$-vacant.
In a potentially misleading  but not uncommon usage,
the terms $p$-dual-path, $p$-dual-cluster, etc.\ are
taken here to mean the corresponding objects on the dual  graph.
The vacant edges of MST on a graph $G$
form a random spanning  tree model, which can  be constructed
as MST on $G^*$ with call numbers $u_{b^*}= 1-u_b$.

The next lemma relates the crossings of
$D(r,R)$ by paths in MST to crossings of
the annulus by curves pieced together from
$p_c$-paths and $p_c$-dual paths. Define a {\em $p_c$-semipath}
to be a (oriented) curve consisting
of a  $p_c$-dual path $\C^+$ and a $p_c$-path $C^-$
such that there is a pair of dual edges $b$ and $b^*$, so that
$b^*$ contains the last vertex of  $\C^+$, and $b$ contains
the first vertex  of $\C^-$ as an endpoint. We allow the special cases
of a $p_c$-path (i.e.\ $\C^+$ is empty) or a $p_c$-dual path
($C^-$ is empty).
We say a $p_c$-semipath traverses an annulus $D(r,R)$, if it
connects a (dual) vertex on one boundary of $D(r,R)$ with
a vertex on the other boundary.  Two semipaths  are {\em disjoint} if no edge
or dual  edge of the one is the same or dual to an edge
or dual edge of the other.

\begin{figure}[htb]
    \begin{center}
    \leavevmode
        \epsfysize=2in
  \epsfbox{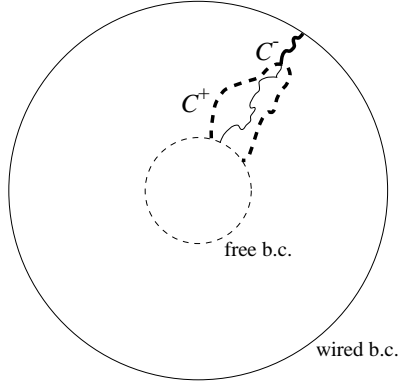}
\caption{\footnotesize
A $p_c$ semipath consists of
a $p_c$-dual path $\C^+$ and a $p_c$-path $\C^-$
joined  at a bond/dual bond pair.
}
\label{fig:semipath}
\end{center}
\end{figure}

\begin{lem} Suppose $\C_1,\dots \C_k$
are disjoint curves in a realization
of MST with mixed (free-wired or wired-free)
boundary conditions on $D(r,R)$ which traverse $D(r,R)$,
where $k\ge 2$. Then the corresponding realization
of Bernoulli percolation contains $k$ disjoint crossings
of the annulus by $p_c$-semipaths.
\label{lem:MST-perc}
\end{lem}

\begin{proof} To be specific, consider the case of free-wired
boundary conditions (the other case is analogous).
Orient the curves $\C_i$ to run from the free boundary at $r$ to
the wired boundary at $R$. If $\C_i$ is a $p_c$-path, then take
$\C_i^-=\C_i$, $\C_i^+=\emptyset$.  For each $i$ such that
$\C_i$ is not a $p_c$-path, let
$b_i$ be the last edge along $\C_i$ with
$u_{b_i}\ge 1/2$. The portion  of $\C_i$ between
$b_i$ and the wired boundary forms a $p_c$-path, which we take to be $\C_i^-$.
By Lemma~\ref{lem:cluster} applied to the dual tree, the two endpoints
of $b_i^*$  are joined to each other by a $p_c$-dual path,
which must pass through $\partial B(r)^*$ because it cannot cross
$\C_i$.  Thus each  of the sectors of the annulus cut out by
the set of $C_i$'s contains  two of these $p_c$-dual paths,
which may well intersect. To obtain a
collection of disjoint $p_c$-semipaths $(\C_i^+,\C_i^-)$ ,
choose $\C_i^+$ to be the $p_c$-dual
path joining $\partial B(r)^*$ to the endpoint of
$b_i^*$ in the sector immediately counterclockwise from $\C_i$.
\end{proof}

One consequence of the lemma is  that for MST on an annulus
of sufficiently large aspect ratio, the probability
of $k$ crossings
decays geometrically in $k$:

\begin{cor} Let $\alpha$ be the exponent defined for critical
Bernoulli  percolation by (\ref{eq:def-alpha-MST}).
For every $s<\alpha/2$,
there exists $m$ large enough so that MST has the geometric-decay
property (\ref{eq:geom-decay}) on annuli $D(r,R=3^{2m}r)$.
\label{cor:geom-decay-MST}
\end{cor}

\begin{proof} We will show that for $r$ and  $R$ as described
in the assertion,
\begin{equation}
\P_\delta\left(\begin{array}{c}
\text{$\Gamma^{F,W}_{r,R} [\Gamma^{W,F}_{r,R}]$ contains $k$ disjoint}\\
\text{traversals of $D(r,R)$}\end{array}\right)
\ \le\ \left( \frac{r}{R} \right)^{sk}\
\quad\text{for all }\ k\ge 2,\ 0\le\delta\le\delta_o(r)\ ,
\label{eq:geom-decay-MST}
\end{equation}
which clearly implies the claim.

Consider the case of free-wired boundary
conditions.  By Lemma~\ref{lem:MST-perc}, there corresponds to
a given collection of at least two tree
crossings $\C_i$ ($i=1,\dots ,k$)
a disjoint collection of $p_c$-semipaths $(\C_i^+,\C_i^-)$,
joined at $b_i$.  Let $n$ be the number
of crossings where either $\C_i$ is a $p_c$-semipath,
or $b_i$ lies in the inner annulus $D^{\text{in}}=D(r,3^mr)$
or else $b_i$ crosses the intermediate boundary at $3^mr$.
Then the semipaths contain $n$ $p_c$-paths traversing
the outer annulus $D^{\text{out}}=D(3^mr,3^{2m}r)$
and $k-n$ $p_c$-dual paths traversing the inner annulus.
We obtain
\begin{eqnarray}
\P_\delta \left(\begin{array}{c}
 \text{$\Gamma_{r,R}^{F,W}$ contains $k$ disjoint}\\
\text{traversals of $D(r,R)$} \end{array}\right)
&\le& \sum_{n\le k}
     \P_\delta\left(\begin{array}{c}
  \text{$D^{\text{in}}$ is traversed by at}\\
\text{least $n$ disjoint $p_c$-paths}
 \end{array}\right)\nonumber\\
&& \qquad
     \times \P_\delta\left(\begin{array}{c}
  \text{$D^{\text{out}}$ is traversed by at least}\\
\text{$k-n$ disjoint $p_c$-dual paths}
 \end{array}\right)\nonumber\\
&\le& (k+1)\,3^{-\alpha m k}\nonumber \\
&\le& \left(\frac{r}{R}\right)^{[(\alpha/2 -1/\log{(R/r)}]\,k}\ ,
\end{eqnarray}
where we have used the independence of events in
$D^{\text{in}}$ and $D^{\text{out}}$ gained from the
decoupling boundary conditions in the first line,
inequality~(\ref{eq:def-alpha-MST}), its dual, and the
telescopic principle  for Bernoulli percolation in
the second line, and $(k+1)\le e^k$ in the last line. The assertion
follows by choosing $R/r=3^{2m}$ sufficiently large.
\end{proof}

The corollary implies that $\gamma(k)\ge\frac{\alpha}{2}  \,k$
for $k\ge 2$.  The relation between the exponents for
MST and Bernoulli percolation can be tightened:

\begin{lem} For MST on $\delta\Z^2$,
the exponents $\gamma(k)$ satisfy
\begin{equation}
\gamma(k)\ \ge\ \min_{n\le k}
\left[\gamma_{B}(n) + \gamma_{B}(k-n)\right]\quad (k\ge 2)\ .
\end{equation}
\label{lem:MST-perc2}
\end{lem}

\begin{proof} Consider, again, MST with free-wired boundary
conditions on $D(r,R)$.  Subdivide $D(r,R)$ into
$M$ annuli $D_j$ of equal aspect ratio $(R/r)^{1/M}$.
By Lemma~\ref{lem:MST-perc}, any collection
of at least two disjoint traversals $\C_i$ of $D(r,R)$ by
$\Gamma_{r,R}^{F,W}$ gives rise to a collection of disjoint traversals
by $p_c$-semipaths $(\C_i^+,\C_i^-)$. Hence each of the annuli $D_j$ is
traversed by a number $n_j$
of $p_c$-paths and at least $k-n_j$ $p_c$-dual paths, with the
possible exception of at most $k$ annuli which
meet one of the special edges $b_i$  (if $b_i$  crosses the
boundary  between $D_j$ and $D_{j+1}$, we discard only $D_j$.)
Let $A_j^-$ (resp. $A_j^+$) denote  the event that $D_j$ is
traversed by $n_j$ disjoint $p_c$-paths (resp., by $k-n_j$ disjoint
$p_c$-dual paths). Then, by the FKG inequalities,
\begin{equation}
\P (A_j^-\cap A_j^+)\ \le\ \P(A_j^-) \,\P(A_j^+)\ .
\label{eq:FKG}
\end{equation}
Using this after  summing over the possible positions
of the $b_i$, and using the independence
of spatially separated events
as in the proof of Corollary~\ref{cor:geom-decay-MST} we obtain
\begin{eqnarray*}
\P_\delta \left(\begin{array}{c}
\text{$\Gamma_{r,R}^{F,W}$ contains $k$ disjoint}\\
\text{traversals of $D(r,R)$}\end{array}\right)
&\le& \P_\delta\left(\begin{array}{c}
\text{$D(r,R)$ is traversed by}\\
\text{at least $k$ $p_c$-semipaths} \end{array}\right)\\
&\le& M^k\, \left(\frac{r}{R}\right)^{(1-k/M)\,
           \min_{n\le k}[\gamma_B(n) +\gamma_B(k-n)]}\ .
\end{eqnarray*}
Choosing $M$ sufficiently large proves the claim.
\end{proof}

\begin{cor} {\rm({\bf H1} for MST)} For all $k\ge 2$,
\begin{equation}
\gamma(k) \ge \ \alpha\,k \quad ,
\label{eq:lin-gamma-MST}
\end{equation}
where $\alpha>0$ is the exponent defined for critical Bernoulli percolation
by (\ref{eq:def-alpha-MST}). In particular,
{\bf H1} holds for MST with
\begin{equation}
\lambda(k) \ \ge \lamst(k)\ \ge \ \phi(k)\ \ge\ \frac{\alpha}{2}(k-1) \   .
\label{eq:lin-lambda-MST}
\end{equation}
\label{cor:H1-MST}
\end{cor}

\begin{proof} Just combine Lemma~\ref{lem:MST-perc2} with
(\ref{eq:H1-MST}),  and with the results
of Lemmas~\ref{lem:lambda-phi} and
\ref{lem:phi-gamma}.
\end{proof}

\masubsect{Euclidean spanning tree}

The proof of {\bf H1} for EST follows the same general strategy as
the proof for MST in the previous subsection. The basic idea
is to relate the tree process to a percolation process,
in this case {\em droplet} percolation (sometimes called
{continuum} or {lily-pad} percolation).  There are a few
additional difficulties, related
with  the lack of self-duality, and the fact that events
in disjoint, but neighboring regions need not be independent.
As a consequence, the definition of {\em disjointness}
for  dual traversals becomes more complicated, and the
relation we establish between crossing events in EST and droplet
percolation is not so tight.  But let us now turn to
the details.

In the introduction, we  defined EST in $\R^2$ as the minimal
spanning subtree of the complete graph on a collection of Poisson points
with density $\delta^{-2}$, with the edge length given by
Euclidean distance. In the droplet percolation model, the random objects
of interest are the connected clusters formed by discs
of a fixed radius  $p\delta$ (where $p$ is a parameter)
centered on the Poisson points.  By construction,
the Poisson process defines a coupling of EST to
droplet percolation with any parameter value $p>0$.

A {$p$-path} is a simple polygonal curve whose straight
line segments join Poisson points with distance
less than $2p\delta$.  A {\em $p$-cluster} is a maximal
set of points that can be joined by
$p$-paths.  As in the case of Bernoulli percolation,
there is a critical value $p_c$ for the parameter.
It follows from the results
of~\cite{Alex-RSW}  (see in particular the proof of Theorem 3.4 and
Corollary 3.5 there) that
for annuli of some fixed aspect ratio $\sigma$,
\begin{equation}
\P_\delta \left( \begin{array}{c}
\text{ $D(r,\sigma r)$ is traversed}\\
\text{by a $p_c$-path} \end{array} \right)\ \le\
\sigma^{-\alpha} \quad (0<\delta\le\delta_o(r))\ ,
\label{eq:def-alpha-EST}
\end{equation}
with some $\alpha>0$. Two $p$-paths or two paths in EST are regarded
as {\em disjoint}, if they share none of their Poisson points.
With this notion of disjointness, a van den
Berg-Kesten inequality
holds for the probability of multiple disjoint $p$-crossings,
and  we obtain as in the case of Bernoulli percolation
the geometric-decay property
\begin{equation}
\P_\delta \left(\begin{array}{c}
\text{$D(r,\sigma r)$ is traversed by}\\
\text{at least $k$ disjoint $p_c$-paths}\end{array}\right)
\ \le\ \sigma^{-\alpha k}\ .
\label{eq:geom-decay-droplet}
\end{equation}
A telescopic argument as in Lemma~\ref{lem:telescopic}
implies that {\bf H1} holds for droplet percolation
in $\R^2$, with exponents  $\gamma_{D}(k) \ge \alpha k$.

One notable difference to Bernoulli percolation is that
droplet percolation is not self-dual. A $p$-dual cluster
is a vacant space inside which a
disc of radius $p\delta$ can be moved without
touching any Poisson points. A $p$-vacant curve is a simple
curve which keeps a distance of at least $p\delta$ to all Poisson points.
The results of \cite{Alex-RSW} imply that
\begin{equation}
\P_\delta \left( \begin{array}{c}
\text{ $D(r,\sigma r)$ is traversed}\\
\text{by a $p_c$-vacant curve} \end{array} \right)\ \le\
\sigma^{-\alpha^*} \quad (0<\delta\le\delta_o(r))\ ,
\label{eq:def-alpha*-EST}
\end{equation}
with some $\alpha^*>0$. (We have chosen $\sigma$ large
enough so that the same $\sigma$ may be used
in (\ref{eq:def-alpha-EST}) and (\ref{eq:def-alpha*-EST}).)  From
this, a geometric-decay property
can be obtained for multiple crossing events ---
if a van den Berg-Kesten inequality is available.
In order to extend the van den Berg-Kesten inequality
from Bernoulli random variables to the present context,
we define a very strict notion of disjointness:
Two $p$-vacant curves are {\em spatially separated}, if
their $p\delta$-neighborhoods are disjoint, i.e., if any pair of
points on the two curves has distance at least $2p\delta$.
Then
\begin{equation}
\P_\delta \left(\begin{array}{c}
\text{$D(r,\sigma r)$ is traversed by at least}\\
\text{$k$ spatially separated $p_c$-vacant curves}\end{array}\right)
\ \le\ \sigma^{-\alpha^* k}\ .
\label{eq:geom-decay-vacant}
\end{equation}
so that ${\bf H1}$ holds for vacant percolation
with exponents $\gamma_{D}^*(k)\ge k\alpha^*$,
whose value may differ from the parameters for the droplet
percolation model itself.

As mentioned in the introduction, EST is automatically a
subgraph of the Poisson-Voronoi graph~\cite{Prep-Sha} with  the natural
Euclidean edge lengths. It can be constructed with the
invasion algorithm of the previous subsection, with any vertex as
the root. An edge of the Poisson-Voronoi graph  will be called 
{\em p-occupied}
if it joins a pair of Poisson points of distance at most $2p\delta$,
and  {\em p-vacant} otherwise. Clearly, Lemma~\ref{lem:cluster}
continues to hold for  $EST$ in place of $MST$, with $\delta\Z^2$ replaced
by the Poisson-Voronoi graph of density $\delta^{-2}$ on $\R^2$,
and Bernoulli percolation replaced by droplet percolation.

For any random spanning tree model on a planar 
graph $G$, we can construct a dual tree model on the dual graph $G^*$,
as explained in the previous subsection.  The dual of a Poisson-Voronoi
graph in $\R^2$ can be represented with the corners of the
Poisson-Voronoi cells as dual vertices,
and the straight line segments of the cell boundaries
as dual edges. A $p$-dual path   
is a simple polygonal curve  consisting of the duals
of {\em p-vacant} edges in $G^*$, i.e., of boundaries of
cells defined by Poisson points that are at least a distance  
$2p\delta$
apart.  (See the discussion of MST for the effect of
free and wired boundaries.) Since
a $p_c$-dual path in $G^*$
is clearly a $p_c$-vacant curve, Lemma~\ref{lem:cluster} holds
also for the dual of EST and vacant percolation
(in place of MST and  and Bernoulli percolation, respectively).

In accordance with the previous
definition, we define a $p_c$-semipath $(\C^+,\C^-)$
in the Poisson-Voronoi graph of density $\delta^{-2}$
to be a (oriented) curve consisting
of a $p_c$-dual path in $G^*$,
and a $p_c$-path $\C^-$ in $G$ such that
the last dual vertex of $\C^+$
lies in the boundary of the cell containing the  first
vertex  of $\C^-$. (We allow the same special cases
as before.) Tightening the previous definition,
we say that two $p_c$-semipaths  are {\em disjoint}
if they share no vertices or dual vertices.
Then Lemma~\ref{lem:MST-perc}
continues to hold for EST in place of MST.

Although a $p_c$-dual path in $G^*$ always defines a 
 $p_c$-vacant curve in the plane, and conversely,
a $p_c$-vacant curve can be deformed to run
along the boundaries of Voronoi cells, the notions of
disjointness (of $p_c$-dual paths in $G^*$)
and of spatial separation (of $p_c$-vacant curves in the plane)
are different, and our proof of
 Corollary~\ref{cor:geom-decay-MST} has to be changed
accordingly:

\begin{cor} Let $\alpha$, $\alpha^*$, and $\sigma$
be the parameters defined for droplet and vacant
percolation in (\ref{eq:def-alpha-EST}) and
(\ref{eq:def-alpha*-EST}). For $s\le \min(\alpha,\alpha^*)/4$,
EST has the geometric-decay property (\ref{eq:geom-decay})
on annuli $D(r,R=\sigma^{2m}r)$ with a sufficiently large integer $m$.
\label{cor:geom-decay-EST}
\end{cor}

\begin{proof} We will show that, for $r$ and $R$ as in the statement,
\begin{equation}
\P_\delta \left(\begin{array}{c}
\text{$\Gamma_{r,R}^{F,W}[\Gamma^{W,F}_{r,R}]$ contains $k$ disjoint} \\
\text{traversals of $D(r,R)$}
\end{array}\right)
\ \le\  \left(\frac{r}{R}\right)^{2s\lfloor k/2 \rfloor }
\quad \text{for all}\ k\ge 2,\ 0<\delta\le\delta_o(r)\  .
\end{equation}

Subdivide $D(r,R)$ into
an inner annulus $D^{\text{in}}=D(r,\sigma^mr)$ and an outer annulus
$D^{\text{out}}=D(\sigma^mr,\sigma^{2m}r)$, and consider
the disjoint semipaths
$(\C_i^+,\C_i^-)$  corresponding to the $k$ traversals of
the annulus by the tree. As in the proof of
Corollary~\ref{cor:geom-decay-MST}, we obtain
$n$ crossings of  $D^{\text{out}}$
by $p_c$-paths $\C_i^-$ in the Poisson-Voronoi graph,
and $k-n$ crossings $D^{\text{in}}$ by $p_c$-dual paths $\C_i^+$.

By definition, each $\C_i^-$  is a $p_c$-path for droplet
percolation, and disjoint $p_c$-semipaths lead to disjoint
$p_c$-paths.  Similarly, each of the paths $\C_i^+$
along the edges in $G^*$
can be parametrized as a curve in the plane that
keeps distance  at least $p_c\delta$ from all Poisson points.
The complication  here is that the $p_c$-vacant
curves $\C_i^+$ need not be spatially separated according to our
definition given above even for disjoint semipaths.
However, by Lemma~\ref{lem:vacant-separated} proved below,
we can use the way the $\C_i^+$ are confined to the sectors
cut out of $D(r,R)$ by the set
of $\C_i$'s, to find at least $\lfloor (k-n)/2\rfloor$ $p_c$-vacant paths
among the $\C_i^+$'s which are spatially separated, except possibly,
for their first and last edges. (As usual, the possibility of
long edges introduces a correction which  is exponentially small
in $\delta^{-2}$.)

The proof is completed by using the
independence of events in $D^{\text{in}}$ and $D^{\text{out}}$
(with the decoupling boundary conditions), and
the geometric decay properties (\ref{eq:geom-decay-droplet})
and (\ref{eq:geom-decay-vacant}) for droplet and vacant percolation.
\end{proof}

\begin{lem}  Let $b=\{x,y\}$ be  an edge of
EST with density $\delta^{-2}$, and let $P$ be a $p_c$-dual path
in $G^*$ (the dual of the corresponding  Poisson-Voronoi graph) 
with $p_c$  the critical parameter value for
droplet percolation.  Assume that no edge of
$P$ is dual to $b$. Then the distance between $b$
and all non-terminal segments of $P$ is at least $p_c\delta/2$.
\label{lem:vacant-separated}
\end{lem}

\begin{figure}[htb]
    \begin{center}
    \leavevmode
        \epsfysize=1in
  \epsfbox{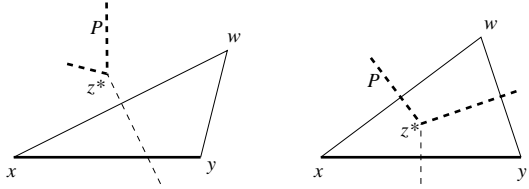}
\caption{\footnotesize
Two possible positions of an edge $b=\{x,y\}$ in the Poisson-Voronoi
graph relative to a $p_c$-dual path $P$ containing the
boundaries of the Voronoi cells of $x$ and $y$. The cells of
$x$ and  $y$ meet the cell of $w$ at $z^*$, which is point 
on $P$ closest to $b$.
}
\label{fig:triangle}
\end{center}
\end{figure}

\begin{proof} The minimal
distance between $b$ and the non-terminal segments of $P$ is
realized
for a pair of points involving either an endpoint of
$b$ or the endpoint of a segment of $P$.  In the first case, we are 
done, since $P$ has distance at least $p_c\delta$ from
any Poisson point, and in particular from  the vertices
$x$ and $y$.  In the second case, the minimal distance is assumed
somewhere between a point on $b$ and a vertex $z^*$ of $G^*$
on $P$.  We need to find a lower bound for the height
$h$ of the triangle  $xyz^*$. Assume, without loss of generality,
that $z^*$ lies on the common boundary of the
Voronoi cells of $x$ and $y$ with the cell of another point $w$ (otherwise,
the tree contains an edge that is closer to $z^*$ than $b$).
In other words, $z^*$ is the  center of the circle
through $x$, $y$, and $w$.

Both $\{x,w\}$ and $\{y,w\}$ have
length at least $2p_c\delta$, because $P$ contains their duals.
Moreover, one of them (say $\{x,w\}$) is longer than $b$, 
because EST contains $b$.  If the triangle $xyw$  has an obtuse angle
at $y$, then $\{x,w\}$ has  length
at least $\sqrt{4(p_c\delta)^2+\ell^2}$ (where $\ell$ is the length of $b$),
so that the distance of $z^*$ to both $x$ and $y$ exceeds half of that value.
Since $z^*$ lies on the perpendicular bisector of $b$,
we see with the Pythagorean theorem that $h\ge p_c\delta$.

If the triangle $xyw$ has acute angles at both $x$ and $y$,
we slide $x$ and $y$ apart in such a way that the line through
$x$ and $y$ and their perpendicular  bisector
are preserved, until the lengths of $\{x,y\}$ and $\{x,w\}$
coincide. While this increases the lengths of all sides
of the triangle $xyw$, it can only decrease $h$,
since the intersection of the Voronoi cells of $x$ and $y$ with
the bisector
of $b$ shrinks. Elementary geometric considerations show that
$h^*$ exceeds $p_c \delta /2$ (see Figure~\ref{fig:triangle}).
\end{proof}

Lemma~\ref{lem:MST-perc2} and Corollary~\ref{cor:H1-MST}
have to be modified as well:

\begin{lem} In the case of EST of density $\delta^{-2}$ on $\R^2$,
the exponents $\gamma(k)$ satisfy
\begin{equation}
\gamma(k)\ \ge\ \min_{n\le k}
\left[\gamma_D(n) + \gamma_D^*(\lfloor (k-n)/2\rfloor )\right]
                 \quad (k\ge 2)\  .
\end{equation}
\label{lem:EST-perc2}
\end{lem}

\begin{cor} {\rm({\bf H1} for EST)} For all $k\ge 2$,
\begin{equation}
\gamma(k) \ge \ \min(\alpha,\alpha^*)\,
           \left\lfloor \frac{k}{2}\right\rfloor \quad ,
\label{eq:lin-gamma-EST}
\end{equation}
with $\alpha, \alpha^*>0$ as defined above. In particular,
{\bf H1} holds for EST with
\begin{equation}
\lambda(k) \ \ge \lamst(k)\ \ge \ \phi(k)\ \ge\ \min(\alpha,\alpha^*)\,
          \left\lceil\frac{k-1}{4}\right\rceil\quad  .
\label{eq:lin-lambda-EST}
\end{equation}
\label{cor:H1-EST}
\end{cor}

\begin{proof}  Combine Lemma~\ref{lem:EST-perc2} with
the general inequalities between the exponents of
Lemmas~\ref{lem:lambda-phi} and \ref{lem:phi-gamma}.
\end{proof}

\remark In Corollary~\ref{cor:H1-EST} and
Lemma~\ref{lem:EST-perc2}, the expression $\lfloor (k-n)/2\rfloor$
can be replaced by $1$ when $k-n=1$.

\masect{Verification of {\bf H2} }
\label{sect:rough}


We shall now verify the roughness criterion.  In contrast
with the
previous section, our  arguments here will rely mostly 
on the tree structure, symmetry, and planarity. In particular,
the result of this section also applies to the
uniform spanning tree on the
Poisson-Voronoi graph.  The main idea is seen
in the following lemma.

\begin{lem}
Let $\Gamma(\omega)$ be a random tree model on $\R^2$, and let $B$ 
be a rectangle in the plane. Suppose that the distribution of
the model is symmetric under a group of
transformations in the plane which is large enough
so that some collection $B_1,\dots ,B_n$ of images of $B$ under
these transformations can be positioned  in such a way that any
collection of $n$ curves $\C_i$ traversing $B_i$
($i=1,\dots , n)$ forms a loop.
Then
\begin{equation}
\P\left(
\begin{array}{c}
        \text{ $B$ is traversed (in the long } \\
        \text{ direction) by a path in $\Gamma$ }
\end{array}
  \right) \ \le \ 1 - \frac{1}{n}  \;  .
\label{eq:h2event}
\end{equation}
\label{lem:h2event}
\end{lem}

\remark If the model has the symmetries of the
square lattice, $B$ can be any sufficiently  long rectangle, and
$B_i$
($i=1,\dots ,4$) are the images of $B$ under rotation by $\pi/2$
about a sufficiently close lattice point. For the hexagonal and
triangular lattice, we would use rotations by $2\pi/3$
in the same way.

\begin{proof}
Since the random tree contains no loops, the probability
that all $B_i$ are traversed simultaneously must vanish.
Thus, with probability one at least one of the $B_i$
fails to be traversed.  By our symmetry assumption,
the probability of failure has to be at least
$1/n$, which proves \eq{eq:h2event}.
\end{proof}

The above observation will now be supplemented by a decoupling 
argument.  

\begin{lem} Let $\Gamma_\delta(\omega)$
be one of the four spanning tree models on $\R^2$  described in the
introduction, with cutoff parameter $\delta$,  and let 
$\{ A_1,\ldots,A_k \}$ be a collection of well separated rectangles 
of common aspect ratio (length/width) $\sigma > 2$.  Then 
\begin{equation}
\limsup_{\delta \to 0} \ \P_\delta\left( \begin{array}{c}
  \text{each $A_j$ is traversed (``lengthwise'') } \\
  \text{by a curve in $\F_{\delta}^{(2)}(\omega)$ }
  \end{array}
\right)\ \le\  \rho^k \;   
\end{equation}
with $\rho = 3/4$.  Furthermore, with some other values of 
$\rho < 1$, and $\sigma < \infty$, the above bound on the probability 
applies for all $\delta < \min_j \ell_j$ (i.e., 
also the full hypothesis {\bf H2} holds). 
\label{lem:H2}
\end{lem} 
\begin{proof} \ 
Let us consider first the case of the spanning trees on  $\delta \Z^2$. 
For each of the $A_j$, we pick a lattice point  $x_j$
outside $A_j$, but as close as possible to the midpoint
of one of the long sides.  Let $\Lambda_i$ be
the disc of radius $\sigma \ell_j$ about $x_j$.
Then $\Lambda_j$ contains $A_j$.  The discs are disjoint
since the separation between
$A_j$ and the other rectangles is larger
than $2\sigma\ell_j$.
Introducing free boundary conditions on the
$\Lambda_j$ will only enhance
the crossing probabilities, while decoupling the crossing events
in disjoint discs.  We next check the assumptions
of Lemma~\ref{lem:h2event}. Clearly, in each
of $\Lambda_j$  the tree processes   
(with free boundary conditions)
is symmetric under rotation by $\pi/4$ about $x_j$.
If $\delta$ is small enough
($ \delta  \le \frac{\sigma-2}{4\sqrt{2}}
\,\min{\ell_j}$ will do), then the images
of $A_j$ under the four rotations by
multiples of $\pi/2$ intersect in such  a way that
any simultaneous crossings would form a loop.
By Lemma~\ref{lem:h2event} the crossing probabilities 
are independenty bounded by $3/4$. This implies both claims for 
the UST and the MST.

An additional consideration is needed for the models 
on the Poisson-Voronoi graph.   One may
take here $\sigma=2$, choose $x_j$ to be the midpoint of
a long side of $A_j$, and let $\Lambda_j$ be the disc of
radius $\sigma\ell_j$ about $x_j$.  
The 
the probability that each $A_j$ is crossed by the
restriction of the tree to $\Lambda_j$ is bounded by $(3/4)^k$, 
by the same argument as above.  However, a small 
correction has to be added to allow for the possibility
of an edge crossing $A_j$ and the boundary of $\Lambda_j$.
As discussed in Section~\ref{sect:crossing}, the probability
of such a long edge can be dominated by 
$B e^{-A(\sigma \ell_j/\delta)^2}$, with suitable constants 
$0 < A,\ B < \infty$.  The claim then easily follows 
also for that case.  
\end{proof}

\masect{Conclusion}
\label{sect:conclusion}


 The scale invariant bounds derived in Sections~\ref{sect:reg}
 and~\ref{sect:rough}  will now be used to prove
 the two Theorems stated in the Introduction.

\masubsect{Tightness, regularity, and roughness}
\label{subsect:limit}

The basic strategy for the proof of Theorem~\ref{thm:main1}
is to apply the regularity and roughness
results for random curves
(Theorems~\ref{thm:reg-curves} and~\ref{thm:rough-curves},
see Section~\ref{sect:criteria})
to the branches of the random trees to obtain the
tightness of the family $\{\mu_\delta^{(2)}\}$, and then use
the structure of the spaces $\Om^{(N)}$ and $\Om$
to obtain tightness of $\{\mu_\delta^{(N)}\}$ and $\mu_\delta$.
The statement about the locality and basic  structure
follows from the positivity of $\lambda(2)$.

\bigskip
\begin{proof_of}{Theorem~\ref{thm:main1}} \

\medskip
{\em Existence of limit points:} We verified that
$\F^{(2)}_\delta$  satisfies the regularity  criterion {\bf H1}
in $\R^2$
for each of the systems of curves along UST, MST, and EST
(Corollaries~\ref{cor:H1-UST}, \ref{cor:H1-MST}, and~\ref{cor:H1-EST},
respectively). By Lemma~\ref{lem:H1-sphere}, the corresponding
bound on crossing probabilities
holds (with the same exponents) also for the 
system on $\sp^2$ with the metric $d(x,y)$ given by
(\ref{eq:def-metric}).
Theorem~\ref{thm:reg-curves} implies that
the family of measures $\mu^{(2)}_\delta$ is tight, and 
that subsequential scaling limits exist for  the system of random 
curves $\F^{(2)}$.
Since for $N>2$ the spaces $\S^{(N)}$, constructed
by patching together spaces $S^{\tau}$, are closed
subspaces of $\left[\S^{(2)}\right]^{2N-3}$ (see the
discussion at the end of Subsection~\ref{sect:graphs}.a),
the family of measures $\mu_\delta^{(N)}$  on $\Om^{(N)}$
is tight also for each $N>2$. (There is nothing to show for $N=1$.)
Tightness of the measures $\mu_\delta$
on the product space $\Om\subset{\sf X}_{N\ge 1} \Om^{(N)}$
now easily follows by an application of Tychonoff's theorem.

\medskip
The tightness described above
guarantees the existence of a sequence $\delta_n \to 0$ for which
the limit $\lim_{n\to \infty}\mu_{\delta_n}(\cdot)$
exists in the sense of weak convergence of measures on
${\sf X}_{N\ge 1} \Om^{(N)}$, as described
by \eq{eq:weakconvergence}.

To see that a limiting configuration typically describes
a single spanning tree  in $\R^2$, we use
that the exponent $\lambda(2)$  is positive
by Corollaries~\ref{cor:H1-UST}, \ref{cor:H1-MST}, 
and~\ref{cor:H1-EST}.
For $r>0$ and $\delta>0$,  
define the random variable $R_{\delta;r}(\omega)$ 
to be the radius of the smallest ball containing
all trees with endpoints in $B(r)$ , and let $R_r(\omega)$
be the corresponding variable in a scaling limit.
Condition {\bf H1} says that
\begin{equation}
\P_\delta \left( \frac{R_{\delta;r}(\omega)}{r}\ \ge\ u\right)\ \le \ 
K(2,s) u^{-(\lambda(2)-s)}\ ,
\label{eq:quasilocal}
\end{equation}
so that $\F_\delta$ is {\em uniformly quasilocal}
in the sense that $R_{\delta;r}$ is stochastically bounded 
as $\delta\to 0$. Moreover, (\ref{eq:quasilocal})
also holds for $R_{r}(\omega)$ for any scaling limit of 
the system.
In particular, $\mu$-almost every limiting configuration
$\F(\omega)$ is quasilocal, and represents a single tree spanning $\R^2$.

\medskip{\em Regularity:}
Theorem~\ref{thm:reg-curves} guarantees furthermore
that for every $\alpha<1/2$,
the curves in the limiting object $\F(\omega)$
can be parametrized, by functions $g(t)$ which
are H\"older continuous (using the metric given by 
(\ref{eq:def-metric})
on $\sp^2$), with exponent $\alpha$
and a random prefactor  whose distribution depends on $\alpha$,
that is,
\begin{equation}
d(g(t), g(t'))\ \le\ K_{\alpha}(\omega)  
|t-t'|^\alpha\ \qquad 0\le t, t'\le 1  \; . 
\label{eq:holder-sphere}
\end{equation}
Rewriting equation~(\ref{eq:holder-sphere}) in terms of
the original metric on $\R^2$, we obtain
\begin{equation}
|g(t)- g(t')| \  \le \ K_{\alpha}(\omega) \,
   (1+|g(t)|^2 + |g(t')|^2)\, |t-t'|^\alpha  \;  . 
\end{equation} 
  
The last conclusion from Theorem~\ref{thm:reg-curves}
is that in $\mu$-almost all configurations of any scaling
limit, all the curves have Hausdorff dimension at 
most $2-\lambda(2)$.

\medskip{\em Roughness:}
Since $\F^{(2)}_{\delta}$
also  satisfies the roughness criterion {\bf H2$^*$}
by Lemma~\ref{lem:H2}, Theorem~\ref{thm:rough-curves} implies
that the limiting measure $\mu^{(2)}$ is
supported on collections containing
only curves  whose Hausdorff dimension is  bounded below by
some $d_{\min}>1$, which depends on the parameters
in {\bf H2$^*$}. In particular, curves in scaling
limits cannot be
parametrized H\"older continuously with any exponent $\alpha>d_{\min}^{-1}$.
This concludes the proof of the convergence, regularity, and
roughness assertions of Theorem~\ref{thm:main1}.
\end{proof_of}

\masubsect{Properties of scaling limits}

The main tool for the proof of Theorem~\ref{thm:main2} is
 the fact that the limiting measure inherits the power
bounds associated with the exponents $\lamst(k)$,
as explained in Theorem~\ref{thm:semicontinuity}.
It is convenient to employ here the following
notion of degree,  which classifies the local behavior
of a collection of trees near a given point $x \in \R^2$.

\begin{df}  The \underline{degree}
of an immersed
 tree at a point $x$ is given by
\begin{equation}
\deg_T(x)\ =\ \sum_{\xi:f(\xi)=x} \deg_\tau(\xi)\ ,
\end{equation}
where $f:\tau\to\R^2$ is a parametrization of $T$ which
is non-constant on every link.
Here $ \deg_\tau(\xi) $ is the branching number of the reference
tree $\tau$ at $\xi$ if $\xi$ is a vertex of $\tau$, and it is
taken to be $2$ if $\xi$ lies on a link of $\tau$.
For a  collection of trees $\F$ immersed in $\R^d$, the degree
at $x$ is
\begin{equation}
\deg_{\F}(x) \ = \ \sup_N \sup_{T\in\F^{(N)}}\ \deg_T(x)\ .
\end{equation}
\label{df:degree}
\end{df}

A more refined notion is that of the  \underline{degree-type} of
$T$ at $x$, which is
the multiset of the summands in the above definition of degree.  The
notions in Definition~\ref{df:point} can be expressed in terms of
degree-type.  For instance, a point of uniqueness
is one whose degree-type has one part for every tree $T$ in $\F$.
A branching point is one with degree-type (for some $T$ in $\F$)
containing a part that is at least 3, and a pinching point is one with
two parts at least 2.


One may note that $\deg_{\F}(x) = 1$ implies that $x$ is a point
of uniqueness.   Such points are also points of continuity,
in the sense seen in the following statement.

\begin{lem} If $\F$ is a closed inclusive collection
of trees representing a single spanning tree in $\R^d$,
and $\eta=\{x_1, \ldots , x_N \} $ is an $N$-tuple
consisting of distinct points of uniqueness, then
$\F$ includes exactly one subtree, denoted $T^{(N)}(\eta)$,
with the set of external vertices given by $\eta$.

Moreover, if the external vertices of a sequence of
trees $\{T_n\}$  in $\S^{(N)}$
satisfy
\begin{equation}
\eta_n \too{n\to\infty} \eta
\label{eq:eta}
\end{equation}
in $\left(\R^d\right)^N$, then
\begin{equation}
T_n \too{n\to\infty}   T^{(N)}(\eta)\
\end{equation}
with respect to the metric on $\S^{(N)}$.
\label{lem:unique}
\end{lem}

\begin{proof}
Assume that $\F$ contains two trees, $T_1$ and $T_2$
with external vertices given by $\eta$. Since $\F$ represents a single
spanning tree, there exists a tree $T$ (parametrized as $f:\tau\to\R^d$)
containing both $T_1$ and $T_2$, with no external vertices
beyond $\eta$.  If  $T_1\not = T_2$, then at least one of the
two trees (say $T_1$) is parametrized under $f$ by a proper subset
$\tau_1$ of $\tau$.  Let $\xi$ be an external vertex
of $\tau$ not contained in $\tau_1$; clearly $x=f(\xi)$ is one
of the points $x_1,\dots ,x_N$ in $\eta$. By assumption,
there exists a point  $\tilde\xi$ in $\tau_1$ with
$f(\tilde\xi)=x$.  Since $\F$ is inclusive, it contains
the  curve obtained by joining $\xi$ to $\tilde\xi$
in $\tau$ and applying $f$.  This is the desired curve
which starts and ends at $x$.

To see the continuity statement, note that the closedness
of $\F$ implies that any limit of a sequence of trees whose
external vertices satisfy the
assumption (\ref{eq:eta}) is certainly contained in $\F$,
and has external vertices $\eta$. The uniqueness result implies the
claim.
\end{proof}

The dimension of the set of the points of
degree $k$ can be estimated in terms of the exponents
$\lamst(k)$.

\begin{lem} Let $\mu(d\F)$ be a probability measure  on $\Om$
describing a  random collection of trees in $\R^d$, and assume
it satisfies the
power-bound (\ref{eq:lambda-bar}), on the probability of
multiple disjoint crossings of annuli, with a family of
exponents $\lamst(k)$.
For each realization $\F$, let
\begin{equation}
A_k(\F)\ =\ \left\{ x\in\R^d\ \mid\ \deg_{\F}(x)\ge
k\right\} \; .
\end{equation}
Then:
\begin{itemize}
\item[i.]
For $\mu$-almost every $\F$ the Hausdorff dimensions of
$A_k(\F)$ satisfy
\begin{equation}
    \dim_{\mathcal H}A_k(\F)  \ \le \ \left( d-\lamst(k) \right)_+  \; ,
\end{equation}
in particular
\begin{equation}
\lamst(k) > 0  \Longrightarrow  A_k(\F)  \
\mbox{is of zero Lebesgue
measure} \; ;
\end{equation}
\item[ii.]
\bea
 \lamst(k)>d  & \Longrightarrow  & A_k(\F)=\emptyset \
 \mbox{for $\mu$-almost every $\F$, i.e.,}   \\
 & & \quad
          \sup_{x\in \R^d} \ \deg_{\F}(x) \ < \ k,\
          \mbox{$\mu$-almost surely} \nonumber  \; .
 \eea
\end{itemize}
\label{lem:dim}
\end{lem}

\begin{proof}   For $R>0$, we denote by $A_{k,R}(\F)$
 the set of all points $x\in \R^{d}$ such
that for all $r \in (0,R)$ the  tree configuration $\F$
exhibits at least $k$
microscopically disjoint traversals of $D(x,r,R)$.
The definition of the degree implies
\begin{equation}
A_k(\F)\ \subset \ \bigcup_{1\ge R>0} A_{k,R}(\F)\ \; ,
\end{equation}
where it suffices to take $R=2^{-j}$, $j=1,2,\ldots$.
By translation invariance (of $\lamst(k)$ and $\dim_{\mathcal H}$),
and the fact that
the Hausdorff dimension of a countable union
of sets of dimension $\le \nu$ does not exceed $\nu$,
it suffices to show that for any given $R<1$
\begin{equation}
 \dim_{\mathcal H}A_{k,R}(\F)\cap [0,1]^d
   \ \le \ d-\lamst(k)  \ .
   \label{eq:dimension}
\end{equation}
Let now $N(k,r,R;\F)$ be the number of  balls of radius
$r$ needed to  cover $A_{k,R}(\F)\cap [0,1]^d$.
Covering the unit square by $\c r^{-d}$ balls of radius $r$,
we see that for any $s<\lamst(k)$,
the expectation  value satisfies
\begin{equation}
 \E\left( N(k,r,R;\F) \right) \le \c(R,s)\,r^{s-d}  \;  .
 \label{eq:hausdorff}
\end{equation}
By Chebysheff's inequality, the random variables
$r^{d-s} N(k,r,R)$ are stochastically bounded
uniformly in $r$. Equation (\ref{eq:dimension})
readily follows.

In case  $d-\lamst(k) < 0$,
the above covering argument implies that the set is almost
surely empty.
\end{proof}

We shall now use the above observations to complete the proof
of the second set of results stated in the introduction.

\medskip 
\newpage 
\begin{proof_of}{Theorem~\ref{thm:main2}} \

{\em Singly connected to infinity:}
Let $\F(\omega)$ be a scaling limit
of one of the three  random tree
models  considered here (UST, MST, or EST).  Note that
if $\F$ was not singly connected to infinity, then,
with positive probability, it would contain two microscopically
disjoint paths traversing annuli $D(r,R)$ with arbitrary large
aspect ratio. This contradicts the strict positivity of
$\lamst(2)$.

\medskip{\em Points of uniqueness and exceptional points:}
Points of degree one are automatically points of uniqueness.
Thus, the claim that Lebesgue-almost all points are
points of uniqueness is implied by the condition
$\lamst(2) > 0$, through  Lemma~\ref{lem:dim}
with $k=2$.  This also shows that the set of
exceptional points has dimension less than two.

To see  that exceptional points are
dense, it is instructive to consider the dual model, which
in two dimensions is also a spanning tree.
Any interior point of a curve in a
scaling limit of the dual tree model is a point of
non-uniqueness for the original spanning tree.
In two dimensions, the exponents $\gamma(k)$ are shared by
the model and its dual for all the models discussed here
(because the graph $G_{r,R}^{F,W}$ is dual to $(G^*)^{W,F}_{r,R}$),
even in the absence of the self-duality exhibited by UST and MST
so that the dual models also
satisfies  the hypothesis {\bf H1} and {\bf H2$^*$}.
That makes the roughness assertion (\ref{eq:dmin}) of
Theorem~\ref{thm:main1} applicable also to the
dual models, and hence almost surely the dimension
of each dual curve is strictly larger than one.
Also, since  a scaling limit of the dual model is
a single spanning tree, the set of interior
points of its curves is clearly dense in $\R^2$.

\medskip
{\em Countable number of  branching points:}
In order to establish that the collection of
branching points  is countable, it suffices to
show that for every $\eps>0$ there are only countably many
points at which branching occurs with three or more branches
extending to a distance greater that $2\eps$. (The collection
of branching points is a countable union of
such sets, with $\eps=2^{-n}$.) We shall refer to such
points as {\em branching points of scale $\eps$}. As a further
 reduction,
we note that  it suffices to prove that in any finite region,
there are typically only finitely many such points. Thus,
the countability is implied by part  (i)  of
the following claim.

\begin{quotation}
\noindent {\bf Claim:} {\em
Let $N_\eps (\F)$ be the number of points of branching
of scale $\eps$, within the unit cell $\Lambda = [0,1]^2$.  Then
\begin{itemize}
\item[i.] $\mu(d\F)$ - almost surely
\begin{equation}
N_\eps (\F) \ < \ \infty \; ,
\end{equation}
        \item[ii.] for each integer $k$ such that
                $\lamst(k) > 2 \ (=d)$
\begin{equation}
        \P\left( N_{\eps} \ge m  \right) \ \le \
                \frac{\c(k)}{\eps^{\lamst(k)} } \
                \left( \frac{k}{m} \right)^{ \frac{\lamst(k) - 2}{2} } \; ,
\label{eq:rational}
\end{equation}
  for all $m\ge k/\eps^{2}$ (where $\P$
is with respect to the measure $\mu$).
\end{itemize}
}
    \end{quotation}

\begin{proof_claim}\
 Part (i) is of course implied by  (ii).  To prove (ii),
let us partition the unit square into
square cells of diameter $r\le \eps$,
with $r$ determined by
\begin{equation}
        m  \ = \ k / r^{2}  \;  .
        \end{equation}
        This choice of $r$ guarantees that
if  $N_\eps (\F) \ge m$ then in at least one of the cells
$\F$ has $k$, or more, branching points of scale $\eps$.
Now, if a given cell contains $k$ such points, then $\F$ includes
a subtree which within this cell has $k$ branching points, with all
branches extending further than $2\eps-r\ge\eps$ from the cell's center.

This implies that the annulus concentric with the cell, with
inner radius $r$ and outer radius $\eps$, is traversed by at least
$k+2$ microscopically disjoint curves.
(This topological fact was employed in a vaguely related
context by Burton and Keane~\cite{BuKe}.)
Adding our bounds for
the probabilities for such events
($\c(k) (r/ \eps)^{\lamst(k)}$ for each cell), we get
\begin{equation}
        \P\left( N_{\eps} \ge m  \right) \ \le \
                \frac{1}{r^{2} } \ \c(k) \
                        \left(\frac{ r}{\eps}\right)^{\lamst(k)} \; ,
 \end{equation}
  which leads directly to \eq{eq:rational}.
\end{proof_claim}

\medskip
{\em Non-random bound on the degree of branching points:}
The absence of branching points of arbitrary high degree
is a direct consequence of $\lamst(k)\to\infty \ (k\to\infty)$
by Lemma~\ref{lem:dim} (ii).
\end{proof_of}

\noindent{\bf Remark: }  We conjecture that the maximal
branching number is actually $k=3$.   From the perspective
of this work this is suggested by the
countability of the branching points, which  may be an
indication that
$\lamst(3) = 2 \ (=d)$.
If $\lamst(k)$   is also strictly monotone
in $k$, then $\lamst(4) > 2 \ (=d)$ and the suggested
statement then follows by Lemma~\ref{lem:dim} (ii).
However, neither of the two steps in this argument
has been proven.  We note that both are consistent with
the exact predictions for UST, viewed as the $Q\to 0$ limit
of critical Potts models~\cite{Nienhuis,DS}.

\startappendix

\vspace{1truecm plus 1cm}
\noindent {\Large\bf Appendix}


\maappendix{Quadratic growth of crossing exponents}
\label{sect:quad}

In Section~\ref{sect:reg} it was established that
the crossing exponents $\gamma(k)$
for UST, MST, and EST, grow at least
linearly with $k$, as $k\to\infty$.
We shall now prove that the
growth is even faster: quadratic in $k$.
Our derivation extends the analysis of
ref.~\cite{AizISC} where a similar statement was proved for
independent percolation in $d=2$ dimensions.
It was also suggested there (but not proved)
that the proper generalization, for
dimensions $d$ where $\gamma(k)$ does not vanish, should
be  $ \gamma \asymp k^{d/(d-1)}$.
The improved argument presented here yields such a
lower bound for all dimensions $d\ge 2$.

\noindent {\bf Remark: }
It has been proposed for a number of related problems in two
dimensions that exponents similar to $\gamma(k)$ are given {\it
exactly\/} by a quadratic polynomial in $k$~\cite{Nienhuis,DS}.  In
particular, the prediction for UST (viewed as the $Q=0$ critical
Potts model) is $(k^2-1)/4$.  It would be of interest to see
mathematical methods capable of resolving such issues.

We start by deriving an upper bound on the exponents,
using reasoning analogous to that found in ref.~\cite{AizISC}.

\begin{lem}
The actual rate of growth of $\gamma(k)$ is not faster than order
$k^{d/(d-1)}$ for UST.  In $d=2$ dimensions,
that applies also to MST and EST.
\end{lem}

\begin{proof}
We will show for each of the models that there exists a constant
$\beta<\infty$ so that for all spherical shells $D(r,R)$
(with $0<r<R$), and every integer $k$,
\begin{equation}
\P_\delta \left(
     \begin{array}{c}
      \text{$\Gamma_{r,R}^{F,W}$ contains $k$ disjoint}\\
      \text{crossings of $D(r,R)$}
     \end{array}\right)
\ \ge\
\left(\frac{r}{R}\right)^{\beta k^{d/(d-1)}}\quad
(\ 0<\delta\le\delta_o(r,R) \ ).
\end{equation}

To prove this, we show that with sufficiently high probability
there are $k$ crossing
paths which occur separately within $k$ disjoint conical sectors.
The sectors may open at an angle of the order
 $\c k^{-1/(d-1)}$
(where the constant depends only on $d$).  To decouple the events,
we separate the different sectors by imposing
the wired  boundary conditions on the intra-sector boundaries.
For UST the lower bound follows now from the statement that with
probability at least $(r/R)^{\beta k^{1/(d-1)}}$ (for some
$\beta<\infty$), a random walk, and hence also LERW,
started at a point at the center of the sector's inner
(reflecting) spherical  boundary ($|x|=r$) will leave the
sector through its outer spherical boundary ($|x|=R$).
The statement can be derived by a number of random walk techniques.
For $d=2$ dimensions a harmonic function argument yields such a decay
with $\beta = 1/2 + o(r/R)$, i.e., $\gamma(k) \le k^2/2$.
The calculation can be adapted to higher dimensions, but instead of
presenting it here let us outline a qualitative argument.

The desired random walk estimate can be
obtained by noting that when $k\gg 1$
the region to be crossed looks like a narrow pencil,
which may be subdivided into a series of $O(k^{1/(d-1)} \log (R/r))$
pairwise overlapping subregions of moderate aspect ratio.
If the random walk makes it to the middle portion of the outer
boundary of one of the subregions, it is near the center of the
next subregion, and with probability bounded away from 0 will make
it to the middle portion of the outer boundary of the next subregion
without hitting the walls.

For MST and EST in $d=2$ dimensions, we relate the claim to a
crossing event in the associated critical Bernoulli and droplet
percolation models.  Cut $D(r,R)$ into $2k$ sectors of equal width.
(In the case of EST the sectors need to be separated by a gap of
width $2\delta$.)  It was proved in \cite{Russo,SW} that the
probability of finding a $p_c$-crossing (or a $p_c$-dual crossing) in
a given sector is bounded below by $(r/R)^{\beta k}$ with some
$\beta>0$.  Suppose that the configuration of Bernoulli or droplet
percolation has $p_c$-crossings and $p_c$-dual crossings in
alternating sectors.  (By independence of the sectors, this event
occurs with probability $(r/R)^{\beta k^2}$.)

We can construct the tree (MST or EST) associated with the (Bernoulli
or droplet) percolation model via the invasion process described in
Subsection~\ref{sect:reg}.b.  If we start the invasion from any point
where the $p_c$-crossing meets the boundary at $r$, then the invasion
will reach the outer wired boundary before crossing either of the flanking
$p_c$-dual crossings.  Therefore the tree contains a traversal for each
of the $k$ $p_c$-crossings, and these must be pairwise disjoint.
\end{proof}

We proceed to derive a matching lower bound on
the growth rate of the exponents.

\begin{thm}
Suppose a random tree model $\Gamma$ in $d$ dimensions satisfies
the free-wired bracketing principle
\begin{equation}
\Gamma^W_\Lambda \backslash \{ \partial \Lambda \} \  \preceq \
\Gamma \cap \Lambda  \ \preceq  \  \Gamma^F_\Lambda  \; ,
\end{equation}
in a form which yields the telescopic principle with a
negligible error, as in Lemma~\ref{lem:telescopic},
and has the geometric decay property, in the form:
\begin{itemize}
\item[] There exist $\sigma > 1$ and $t> 0$,
such that the random variable $M(r, \sigma; \omega)$ representing
the number of disjoint crossings of the spherical
shell $D(r, R=\sigma r)$ with free-wired boundary conditions, has a finite
moment generating function:
\begin{equation}
E_\delta\left( e^{t \ M(r, \sigma ; \omega)} \right) \ \le
e^{g(\sigma,t)}\quad (0<\delta\le\delta_o(r))
\label{eq:exp}
\end{equation}
with some $g(\sigma,t) < \infty$.
\end{itemize}
Then there exists $\beta > 0$ such that for $R/r$  sufficiently large
\begin{equation}
\P_{\delta}\left( \begin{array}{c}
       \text{$\Gamma^{F,W}_{r,R}$ contains more} \\
\text{than $k$ crossings of $D(r,R)$}
\end{array} \right) \ \le \
K(k,\beta)\,\left( \frac{r}{R} \right)^{\beta\, k^{d/(d-1)}  } \quad
(0<\delta\le\delta_o(r,\beta, k))\ .
\label{eq:quadratic}
\end{equation}
\label{thm:quadratic}
\end{thm}

\noindent {\bf Remarks:}
i) \ Elementary considerations show that the condition (\ref{eq:exp}) is
implied by the geometric-decay hypothesis~(\ref{eq:geom-decay}), which
was derived (for $d=2$) in Section~\ref{sect:crossing}.

ii)i)  It ought to be clear from the proof that the argument
can be extended to other systems, in particular to independent
percolation models and, more generally, to the
Fortuin-Kasteleyn random-cluster models with
$Q \ge 0$ (of course the theorem stated here will be of
interest only for critical states).  For those systems
$M(r, \sigma; \omega)$ will refer to the maximal number
of crossings which can be realized disjointly in the
configuration $\omega$.
The main adjustment needed in the analysis is to replace
the free-wired bracketing principle
by a suitable decoupling boundary condition
which increases the state.  For $0\le Q \le 1$ that is provided
by the free b.c., whereas for $ 1 \le Q $ that role is played
by the wired b.c.  Correspondingly, the assumption made in the
theorem should in each case refer to the statistics of
the variable $M$ under the corresponding b.c.

\begin{proof}  For a given $k$, let us subdivide the spherical
shell $D(r,R)$ into concentric subshells with a common
aspect ratio:
\begin{equation}
D_n \ = \ D(r e^{(n-1) \alpha }, r e^{n \alpha }) \;, \quad
\mbox{ with $\alpha \ = \ b^{-1}k^{-1/(d-1)}$ } \; ,
\end{equation}
where $b>0$ is a parameter whose value will be specified below.
By the telescopic principle,
the probability of $k$ disjoint traversals
of $D(r,R)$ is dominated (up to a negligible error)  by
the product of probabilities of such traversals of the
$ \left\lfloor  \log(R/r)/\alpha \right\rfloor$
subshells, $D_n$, each taken with the decoupling free-wired
boundary conditions.
Thus, as is explained at the end of this proof,  it suffices
to establish the following bound:

\begin{quotation}
\noindent {\bf Claim:}
There are constants $m>0$ and $a(b)\ge 0$, where  $a(b)$ is
strictly positive for small enough $b$,
so that with the above choice of $\alpha$
\begin{equation}
\P_{\delta}\left( \begin{array}{c}
       \text{$\Gamma^{F,W}_{D_1}$ contains more than $k$ disjoint} \\
\text{ traversals   of $D_1=D(r, r e^{\alpha }) $}
\end{array} \right) \ \le \  m e^{- a(b)k}\quad ,
\label{eq:quad-claim}
\end{equation}
for all  $k\ge k_o(b,\sigma,d)$ and $0<\delta<\delta_o(r,k)$.
\end{quotation}

\begin{proof_claim} We
employ a covering of the sphere of radius
$\tilde r = r e^{\alpha/2} $
by balls of radius $ r_o = r \alpha / (2 \sigma ) $
 (see Figure~\ref{fig:beads}),
where $\sigma$ is large enough so that the geometric decay
property~\ref{eq:exp} holds.
Note that
even when the balls are expanded concentrically by the factor
$\sigma$, they do not reach outside $D_1$.
(This can be seen using $1 < e^x - x$ and
$e^x + x < e^{2 x}$ (for $x>0$) with $x=\alpha/2$.)
Thus, each  path crossing $D(r, r e^{\alpha} )$
produces a crossing from the surface of at least one ball
in the cover
to a sphere concentric with it, of radius $r_o \sigma$.
We shall estimate the probability that there are altogether at least
$k$ (or more) such traversals.

By Lemma~\ref{lem:cover} proved below (with $c=r_o/\tilde r$),
there exists a covering of the $\tilde r$-sphere by balls of radius $r_o$,
\[
\tilde r\ S^{d-1}\  \subset\ \bigcup_{x\in A} B(x;r_o)\ ,
\]
which can be partitioned into $A=\cup_{i=1}^{m} A_i$,
in such a way that
\[
B(x;\sigma r_o) \cap B(y;\sigma r_o)\ =\ \emptyset \quad\text{whenever}\
 x,y \in A_i, x\not = y \ .
\]
The important fact is that the partition can be chosen so that
$m=m(\sigma)$ depends only on $\sigma$ and the
dimension (and not on $r_o$ or $\tilde r$). The maximum number of
balls in any of the $A_i$'s is bounded by
\begin{equation}
\max_i\#A_i\ \le \ a_o\,\left(\frac{\tilde r}{r_o}\right)^{d-1}
\ \le\ a_o\,(4b\sigma)^{d-1}\, k\quad
\text{for $k\ge k_o(b,\sigma,\delta)$}
\label{eq:bound-cover}
\end{equation}
provided $k_o$ is large enough so that $\tilde r/r\leq 2$, i.e.,
\ so that
$\sigma r_o/r \le \log{2}$.

\begin{figure}[htb]
    \begin{center}
    \leavevmode
        \epsfysize=2in
  \epsfbox{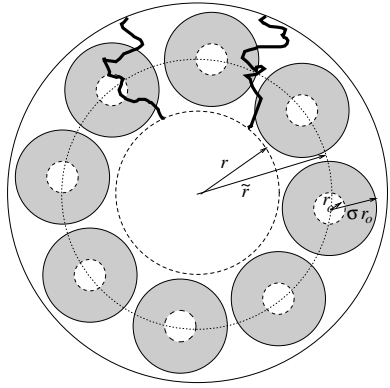}
\caption{\footnotesize
Placement of a disjoint family of small shells
$D(x;r_o, \sigma r_o)$ within a large
shell $D(r,R)$. Of the  two depicted crossings
of the large shell, the left one gives rise to a crossing
of a little shell. In this picture, we have chosen the aspect
ratio $\sigma=3$. Four families of eight disjoint little  shells
each are needed to capture all crossings of the big shells.
}
\label{fig:beads}
\end{center}
\end{figure}

In each configuration, let $M_{i,j}$ be the number of lines
touching the $j$th ball $B(x_j; r_o)$ in $A_i$ (see Figure~\ref{fig:beads}),
and let
$M_i = \sum_j M_{i,j}$.
If there are $k$ disjoint traversals of $D_1$, then
at least one of $M_i$ exceeds $k/m$.  Thus:
\begin{eqnarray}
\P_{\delta}\left( \begin{array}{c}
       \text{more than $k$ crossings} \\
\text{  of $D(r, r e^{\alpha} ) $}
\end{array} \right)
& \le & \sum_{i=1}^m \P_{\delta}\left(  M_i \ge \frac{k}{m} \right)
\nonumber \\
& \le & \sum_{i=1}^m e^{-t k / m}\ \E_{\delta}\left(  e^{t\sum_j M_{i,j}}
\right) \; , \label{eq:Chernoff}
\end{eqnarray}
using in the last step Chebysheff's inequality.
By the free-wired bracketing principle, each of the variables $M_{i,j}$ is
stochastically dominated by the corresponding crossing numbers
$M^{FW}_{i,j}$
of $\Gamma^{F,W}_{r_o,\sigma r_o}$.  We get, for each $i=1,\ldots, m$:
\begin{equation}
\E_{\delta}\left(  e^{t\sum_j M_{i,j}} \right) \ \le\
\E_{\delta}\left(  e^{t\sum_j M^{FW}_{i,j}} \right) \ =\
\prod_j \E_\delta\left( e^{t M^{FW}_{i,j}}   \right) \ \le\
e^{a_o\,(4b \sigma)^{d-1}\, k\,g(\sigma,t) } \; ,
\label{eq:cross-small}
\end{equation}
where we used first the independence of events in disjoint shells
due to the decoupling boundary conditions, and then
the geometric-decay assumption \eq{eq:exp} and
the bound (\ref{eq:bound-cover}) on the number of balls in the cover.
Substituting inequality (\ref{eq:cross-small})
into (\ref{eq:Chernoff}) yields
\begin{equation}
\P_{\delta}\left( \begin{array}{c}
       \text{more than $k$ crossings} \\
\text{  of $D(r, r e^{\alpha} ) $}
\end{array} \right)  \ \le \  m \  e^{- [t / m -
a_o\,(4b \sigma)^{d-1}\,g(\sigma,t)] \cdot k} \; .
\end{equation}
The claim follows now by choosing $\sigma$ and $t$ so that
$g(\sigma,t)$ is finite, and adjusting the parameter $b$, making it
small enough so that
\begin{equation}
a(b) = t / m(\sigma) - a_o\,(4b \sigma)^{d-1}\, g(\sigma,t)  > 0  \; .
\label{eq:diff-exp}
\end{equation}
To make the best out of this argument, one should
optimize in $b$, $\sigma$, and $t$, maximizing $b\times a(b)$.
\end{proof_claim}

The calculation yielding the assertion \eq{eq:quadratic} from the claim
is (for $b$ small enough but independent of $k$)
\begin{eqnarray*}
\P_{\delta}\left( \begin{array}{c}
       \text{more than $k$ disjoint} \\ \text{traversals in
$\Gamma^{F,W}_{r,R}$}
    \end{array}\right)
&\le&  e^{[-a(b) k + \log m] \, \lfloor b\, k^{1/(d-1)} \log(R/r)\rfloor} \\
&\le& e^{a(b) k - \log m} \left(\frac{r}{R}\right)^{b\times a(b)\times
k^{d/(d-1)} - b \times \log m \times k^{1/(d-1)}}.
\end{eqnarray*}
The floors give rise to the prefactor $K(k,\beta)$, which grows
exponentially in $k$.
\end{proof}

Let us remark that the bound (\ref{eq:Chernoff}) makes use
of a standard method for large deviations estimates,
known as Chernoff's inequality (see e.g.~\cite{Chernoff}).
For completeness, following is the covering lemma
used in the analysis.

\begin{lem} {\rm (Covering lemma)} Let $c<1$ and $\sigma>0$.
The unit sphere can be covered with balls of radius $c$
(indexed by a finite set $A$ of centers)
\begin{equation}
S^{d-1}\ \subset \ \bigcup_{x\in A} B(x,c)\ ,
\end{equation}
which can be partitioned into $m$ subcollections
$ A = \bigcup_{i=1}^m A_i$ satisfying
\begin{equation}
B(x,\sigma c) \cap B(y,\sigma c)\ = \
  \emptyset\quad\text{if $x,y\in A_i$, $x\not = y$}\ .
\end{equation}
Here, $m$ depends on $\sigma$ and the  dimension, and
the number of balls needed for the covering  is bounded by
\begin{equation}
\#A \ \le a_o\,c^{1-d} \ ,
\end{equation}
where $a_o$ depends only on the dimension.
\label{lem:cover}
\end{lem}

\begin{proof} In two dimensions, one reasonable choice for $A$ is
a set of evenly spaced points  on the unit circle. In higher
dimensions, take $A$ to be the set of points in
$c\,d^{-1/2}\,\Z^d$ that
are at most distance $c/2$ from the unit sphere.  Every point of
the unit sphere is within distance $c/2$ of such a point.
To bound $\#A$, consider the spherical shell
of inner radius $1-c$ and outer radius $1+c$.  This shell contains
all cubes with side length $c\, d^{-1/2}$ centered about some
point in $A$.  Its volume is bounded above by
$2^d\,\omega_d\,c$ (where $\omega_d$ is the volume of
the unit ball in $\R^d$), so the shell can contain at most
$2^d\,\omega_d\,d^{d/2}\,c^{1-d}$ cubes.  This proves the claim on $\#A$.

By a  similar argument
we conclude that the number of lattice points in a ball of radius
$2\sigma c$ is bounded above by a number $m$ which depends
only on $\sigma$ and the dimension.
We partition $A$ into subsets
$A_1,\dots ,A_m$ so that any two points in $A_i$
have distance at least $s=2\sigma c$, by induction on $m$.
If $m=1$, that is, if the distance between any two points
in $A$ is at least $s$, choose $A_1=A$.  If $m>1$,
take any point and put it in $A_m$ --- this may
make some of the other points ineligible
for placement in $A_m$. Continue in any fashion until all the
points are either in $A_m$, or else ineligible.
Each ineligible point has distance less than $s$ to some
point in $A_m$, so there are at most
$m-1$ other ineligible points at a distance of less than $s$.
Applying the inductive assumption completes the proof.
\end{proof}



\noindent {\large \bf Acknowledgments\/} \\
We thank I. Benjamini, B. Duplantier,
G. Lawler,  R. Lyons, Y. Peres,
and O. Schramm, for stimulating discussions of topics related to this study.
M.A. wishes to thank for the gracious hospitality
accorded him at the School of Physics and Astronomy,
Tel Aviv University, and at the Forschungsinstitut
f\"ur Mathematik, ETH-Zurich,  where some of his work was done.
The work was supported in part by the NSF Grants PHY-9512729 (M.A.),
MPS-9500868 and MPS-9803267 (C.M.N.), DMS-9626198 (A.B.),
and an NSF postdoctoral fellowship (D.B.W.).

\addcontentsline{toc}{section}{References}

\end{document}